\documentclass{elsarticle}
\usepackage{amsmath,amssymb,amsfonts,xcolor}
\journal{Journal of Differential Equations}
\newtheorem{theorem}{Theorem}
\newtheorem{corollary}{Corollary}
\newtheorem{definition}{Definition}
\newtheorem{example}{Example}
\newtheorem{lemma}{Lemma}
\newtheorem{proposition}{Proposition}
\newtheorem{remark}{Remark}
\newproof{proof}{Proof}
\usepackage[lastexercise]{exercise} 
\def\N{\mathbb{N}}

\def\hat{\widehat}
\def\Hat{\widehat}
\def\tilde{\widetilde}
\def\Tilde{\widetilde}
\def\Bar{\overline}
\def\la{\left\langle}
\def\ra{\right\rangle}
\def\ve{\varepsilon}

\def\R{\mathbb{R}}
\def\N{\mathbb{N}}

\def\gph{\mbox{\rm gph}\,}
\def\rep{\mbox{\rm rep}\,}
\def\epi{\mbox{\rm epi}\,}

\def\O{\Omega}

\def\vph{\varphi}
\def\emp{\emptyset}

\def\oR{\Bar{\R}}
\def\lm{\lambda}

\def\gg{\gamma}

\def\dd{\delta}

\def\al{\alpha}

\def\Th{\Theta}

\def\ob{\bar b}

\def\ox{\bar{x}}
\def\oy{\bar{y}}
\def\oz{\bar{z}}
\def\ov{\bar{v}}
\def\ou{\bar{u}}

\def\ow{\bar{w}}
\def\ph{\varphi}

\def\tto{\rightrightarrows}
\def\Tilde{\widetilde}
\def\tilde{\widetilde}

\def\({\left(}
\def\){\right)}
\def\[{\left[}
\def\]{\right]}
\def\n{\left \|}
\def\en{\right \|}

\def\ob{\bar b}

\def\ox{\bar{x}}
\def\oy{\bar{y}}
\def\oz{\bar{z}}
\def\ov{\bar{v}}
\def\ou{\bar{u}}

\def\gph{\hbox{}}

\def\gg{\gamma}

\def\sr{\Longrightarrow }
\def\tto{\rightrightarrows}

\def\hat{\widehat}
\def\Hat{\widehat}
\def\tilde{\widetilde}
\def\Tilde{\widetilde}
\def\Bar{\overline}
\def\la{\langle}
\def\ra{\rangle}
\def\ve{\varepsilon}
\def\bb{\beta}

\def\R{\mathbb{R}}
\def\N{\mathbb{N}}

\def\gph{\mbox{\rm gph}\,}
\def\rep{\mbox{\rm rep}\,}
\def\epi{\mbox{\rm epi}\,}

\def\O{\Omega}

\def\vph{\varphi}
\def\emp{\emptyset}

\def\oR{\Bar{\R}}
\def\lm{\lambda}

\def\gg{\gamma}

\def\dd{\delta}

\def\al{\alpha}

\def\Th{\Theta}

\def\N{I\!\!N}
\def\th{\theta}

\def\disp{\displaystyle}

\def\eq{\begin{equation}}
\def\eeq{\end{equation}}
\begin{document}
\begin{frontmatter}
\author[]{Ren\'e Henrion}
\ead{henrion@wias-berlin.de}
\address{Weierstrass Institute for Applied Analysis and Stochastics, 10117
Berlin, Germany} 
\author[]{Abderrahim Jourani}
\ead{jourani@u-bourgogne.fr}
\address{Institut de Math\'{e}matiques de Bourgogne, {UMR 5584 CNRS, Universit\'{e}
de Bourgogne Franche-Comt\'e, F-21000 Dijon}, France} 
\author[]{Boris S. Mordukhovich}
\ead{aa1086@wayne.edu}
\address{Department of Mathematics, Wayne State University, Detroit, Michigan
48202, USA}

\title{Controlled polyhedral sweeping processes: existence, stability, and
optimality conditions}

\begin{abstract}
\noindent 
This paper is mainly devoted to the study of controlled sweeping
processes with polyhedral moving sets in Hilbert spaces. Based on a detailed analysis of truncated Hausdorff distances between moving polyhedra, we derive new existence and uniqueness theorems for sweeping trajectories corresponding to various classes of control functions acting in moving sets. Then we establish quantitative stability results, which provide efficient estimates on the sweeping trajectory dependence on controls and initial values. Our final topic, accomplished in finite-dimensional state spaces, is deriving new necessary optimality and suboptimality conditions for sweeping control systems with endpoint constrains by using constructive discrete approximations.
\end{abstract}

\end{frontmatter}

\section{Introduction and Problem Formulation}\vspace*{-0.05in}

\noindent In this paper we consider a family of sweeping processes with controlled polyhedral moving sets defined on a Hilbert space $\mathcal{H}$. To describe this family, fix some $x_{0}\in \mathcal{H}$ and, for arbitrary control functions $\left(u,b\right) :\left[ 0,T\right] \rightarrow \mathcal{H}^{m}\times \mathbb{R}^{m}$ satisfying $x_{0}\in C_{\left( u,b\right) }(0)$, define the \emph{moving polyhedral set}
\begin{equation}
C_{\left( u,b\right) }(t):=\left\{ x\in \mathcal{H}|\left\langle
u_{i}(t),x\right\rangle \leq b_{i}(t)\quad \left( \,i=1,\ldots ,m\right)
\right\} \quad \left( t\in \left[ 0,T\right] \right).\label{movpoly}
\end{equation}
This induces the \emph{controlled sweeping process} $\left( \mathcal{S}_{\left( u,b\right) }\right)$ given by
\begin{equation}
-\dot{x}(t)\in
N_{C_{\left( u,b\right) }(t)}\left( x(t)\right) \,\, \mathrm{a.e.\;}\;t\in %
\left[ 0,T\right] ,\,\, x(0)=x_{0}\in C_{\left( u,b\right) }(0),
\label{sweeping}
\end{equation}
where $N_{C}(x)$ stands for the classical normal cone of convex analysis
defined as
\begin{equation}
N_{C}(x):=\{v\in \mathcal{H}\;\big|\;\left\langle v,y-x\right\rangle
\leq 0\}\;\mbox{ if }\;x\in C\;\mbox{ and }\;N_{C}(x):=\emptyset \;%
\mbox{
else}.\label{nc}
\end{equation}
We emphasize that the differential inclusion in (\ref{sweeping}) comes along
with the hidden pointwise \emph{state constraints} $x(t)\in C_{\left( u,b\right) }(t)$ for all $t\in \left[ 0,T\right]$, because otherwise the normal cone is empty by definition.

\emph{Uncontrolled} sweeping processes were introduced and initially studied by Moreau  \cite{JM1,JM2,JM3} and then were extensively developed in the literature, where the main attention was paid to the existence and uniqueness of solutions and various applications; see, e.g., \cite{Adly,Bro,bt,JV,kunzemonteiro,maka} with their references.\vspace*{0.02in}

\emph{Existence and uniqueness} of \emph{class-preserving} solutions $x_{(u,b)}$ to the sweeping dynamics \eqref{sweeping} generated by \emph{control} functions $(u,b)$ in \eqref{movpoly} from various classes in Hilbert spaces is the \emph{first topic} of our paper. Note that the standard approach to this issue (see, e.g., \cite{kunzemonteiro}) consists of checking the Hausdorff Lipschitz continuity of the moving set \eqref{movpoly}. However, this does not make much sense when the moving set is an unbounded polyhedron. The $W^{1,2}$-preserving existence and uniqueness results for moving polyhedra were obtained by Tolstonogov \cite{tolsto,tolsto2,tolsto3} and more recently in \cite{mor1} under certain qualification conditions in Hilbert and finite-dimensional settings; see more discussions in Section~\ref{exsol}. Here we develop a novel approach involving the \emph{truncation} of polyhedra and deriving refined \emph{error bounds}. This allows us obtain new class-preserving results, which shows that Lipschitz continuous (resp.\ absolutely continuous) controls in \eqref{movpoly} uniquely generate Lipschitz continuous (resp.\ absolutely continuous) trajectories of \eqref{sweeping} under an explicit and easily formulated \emph{uniform Slater condition} for moving control polyhedra in separable Hilbert spaces.\vspace*{0.02in} 

The {\em second topic} of our study addresses \emph{quantitative stability} issues on the \emph{H\"olderian} dependence of solutions to \eqref{sweeping} on the corresponding perturbations of controls   $(u,b)$ in moving sets as well as the initial value $x_0$ in separable Hilbert spaces. To the best of our knowledge, such questions have never been posted for the sweeping processes formulated in \eqref{movpoly} and \eqref{sweeping}. Based on the aforementioned truncation techniques and error bounds, we establish efficient results in this direction in the $W^{1,1}$ 
control-trajectory framework.\vspace*{0.02in}

The \emph{third topic} we investigate here concerns an \emph{optimal control} problem for the sweeping process in \eqref{movpoly} and \eqref{sweeping} under the additional pointwise equality constraint on the \emph{$u$-component of controls} and \emph{geometric endpoint constraint} $x_{(u,b)}\in\O$ on trajectories. Optimal control theory for sweeping processes, with addressing the main issue of deriving necessary optimality conditions, has been started rather recently in \cite{chhm} and then has been extensively developed in subsequent publications (see, e.g., \cite{ao,bk,cm,mor1,cg,colmor,mor2,cmn,pfs,zeidan} and the references therein), which did not concern however systems with endpoint constraints. Problems of sweeping optimal control, that are governed by discontinuous differential inclusions with intrinsic pointwise and irregular state constraints, constitute one of the most challenging class in modern control theory. We develop here the {\em method of discrete approximation}, which allows us to constructively approximate the constrained control sweeping process under consideration by discrete-time sweeping systems with perturbed  endpoint constraints so that feasible and optimal solutions to discrete approximations \emph{strongly converge} to the designated feasible and locally optimal solutions of the original problem under the {\em uniform Slater condition} introduced above. Employing then advanced tools of first-order and second-order variational analysis and generalized differentiation, we derive new \emph{necessary optimality conditions} for discrete approximations that gives us efficient \emph{suboptimality conditions} for a general class of local minimizers in the original problem of sweeping optimal control.\vspace{0.02in}

The rest of the paper is organized as follow. Section~\ref{sec:trunc} presents major technical developments on the truncation and error bounds, which are of their own interest while being widely used in deriving the main results of the paper. Section~\ref{exsol} is devoted to establishing the class-preserving existence and uniqueness theorems for the controlled sweeping process. Section~\ref{quantstab} addresses stability issues for sweeping trajectories under control and initial value perturbations. In Section~\ref{discapp} we formulate an optimal control problems for the sweeping process (${\cal S}_{(u,b)}$) with endpoint constraint and construct its well-posed discrete approximations with establishing the $W^{1,2}$-strong convergence of feasible and optimal solutions. The final Section~\ref{sec:optim-disc} provides necessary
optimality and suboptimality conditions for such control problems via advanced tools of generalized differentiation.\vspace*{-0.15in}

\section{Error bounds and truncation of moving sets}\label{sec:trunc}\vspace*{-0.05in}

\noindent This section plays a crucial role in describing and justifying our strategy to derive existence and stability results for sweeping processes with controlled polyhedra in both finite-dimensional and infinite-dimensional settings. The conventional by now theory of sweeping processes establishes the existence of Lipschitz continuous solutions of the sweeping dynamics via the Hausdorff Lipschitz continuity of moving sets; see, e.g., Theorem~2 in \cite{kunzemonteiro} and its proof. Unfortunately, this approach does not work for the case of unbounded moving polyhedra. For instance, in the case in moving \emph{halfspaces}, i.e., for $m=1$ in \eqref{movpoly}, the Hausdorff distance is either zero (if the two halfspaces coincide), or infinity otherwise. Hence the only ``moving" halfspaces satisfying Hausdorff Lipschitz continuity are constant in time, which clearly does not offer any freedom for controlling the process. However, when \emph{truncating} the moving polyhedron with a ball, the Hausdorff Lipschitz continuity may well be achieved. This suggests the following \emph{strategy}, which will be implemented in the paper. \emph{First} we intend to show that Lipschitzian controls lead us to \emph{bounded} continuous solutions of the sweeping process and that the moving polyhedron \emph{truncated} with a ball sufficiently large to contain this solution is Hausdorff Lipschitz, which hence verifies the actual Lipschitz continuity of the solution. The \emph{second step} of our approach is to establish an appropriate \emph{error bound} for the truncation moving polyhedra.

For the reader's convenience, we split this section into several subsections and present numerical examples providing the driving forces for our approach.
\vspace*{-0.07in}

\subsection{\bf Hausdorff Lipschitz continuity of truncated moving polyhedra}

\noindent As discussed above, it is generally hopeless to ensure a
Hausdorff Lipschitz estimate for moving polyhedra \eqref{movpoly} in the form
\begin{equation}
d_{H}\left( C_{\left( u,b\right) }(s),C_{\left( u,b\right) }(t)\right) \leq
\widehat{L}\left\vert s-t\right\vert \quad \forall s,t\in \left[ 0,T\right].
\label{hausorig}
\end{equation}
Our efforts are now paid to establish a \emph{truncated estimate} of type 
\begin{equation}
d_{H}\left( C_{\left( u,b\right) }^{r}(s),C_{\left( u,b\right)
}^{r}(t)\right) \leq \widehat{L}\left\vert s-t\right\vert \quad \forall
s,t\in \left[ 0,T\right],  \label{haustrunc}
\end{equation}
where $r\geq 0$ is appropriately given, and where $C^{r}:=C\cap \mathbb{B}
\left( 0,r\right) $. To accomplish this, we proceed in following two steps. Our \emph{first step} is to derive the \emph{weakened Hausdorff estimate} given by
\begin{equation}
d\left( x,C_{\left( u,b\right) }(t)\right) \leq L\left( \left\Vert
x\right\Vert \right) \left\vert s-t\right\vert \quad \forall s,t\in \left[
0,T\right] \,\,\forall x\in C_{\left( u,b\right) }(s)  \label{weakhaus0}
\end{equation}
with some monotonically increasing function $L(\cdot)$. Estimate \eqref{weakhaus0} clearly yields
\begin{equation}
d\left( x,C_{\left( u,b\right) }(t)\right) \leq \widehat{L}\left\vert
s-t\right\vert \quad \forall s,t\in \left[ 0,T\right] \,\,\forall x\in
C_{\left( u,b\right) }^{r}(s) \label{weakhaus}
\end{equation}
with $\widehat{L}:=L\left( r\right)$. In the \emph{second step} we prove the general estimate
\begin{equation}
d\left( x,C_{\left( u,b\right) }^{r}(t)\right) \leq 3d\left( x,C_{\left(
u,b\right) }(t)\right) \quad \forall t\in \left[ 0,T\right] \,\,\forall x\in
\mathbb{B}\left( 0,r\right)  \label{truncest}
\end{equation}
for all $r$ sufficiently large. Combining the latter with \eqref{weakhaus} will ensure the desired truncated estimate (\ref{haustrunc}). Details follow.\vspace*{-0.05in}

\subsubsection{\bf Limitations of Hoffman's error
bound} 

\noindent The first idea, which comes to our mind for proving (\ref
{weakhaus0}), is the use of the classical {\em Hoffman's error bound}; see, e.g., \cite[Theorem 2.200]{bonnshap}. It guarantees  in our setting that, for each $t\in \left[ 0,T\right] $, there exists some $\widetilde{L}\left( t\right):=L(t,u(t),b(t))$ ensuring the distance estimate
\begin{equation}
d\left( x,C_{\left( u,b\right) }(t)\right) \leq \widetilde{L}\left( t\right)
\max_{i=1,\ldots ,m}\left[ \left\langle u_{i}(t),x\right\rangle -b_{i}(t)%
\right] _{+}\quad \forall x\in \mathcal{H}  \label{hoffman}
\end{equation}%
provided that $C_{\left( u,b\right) }(t)\neq \emp $. In particular, for
$x\in C_{\left( u,b\right) }(s)$ it follows from $\left\langle
u_{i}(s),x\right\rangle \leq b_{i}(s)$ for $i=1,\ldots,m$, that
\begin{eqnarray}
&&\left[ \left\langle u_{i}(t),x\right\rangle -b_{i}(t)\right] _{+} \label{pluest}\\
&=&\left[\left\langle u_{i}(t),x\right\rangle -\left\langle u_{i}(s),x\right\rangle
+\left\langle u_{i}(s),x\right\rangle -b_{i}(s)+b_{i}(s)-b_{i}(t)\right] _{+}
\notag \\
&\leq&\left[ \left\langle u_{i}(t),x\right\rangle -\left\langle
u_{i}(s),x\right\rangle +b_{i}(s)-b_{i}(t)\right] _{+}\notag \\
&\leq&\left\Vert u_{i}(t)-u_{i}(s)\right\Vert \left\Vert x\right\Vert
+\left\vert b_{i}(s)-b_{i}(t)\right\vert \quad \forall i=1,\ldots ,m.\notag
\end{eqnarray}
When $\left( u,b\right)$ is Lipschitz continuous, this combines with the previous estimate to give us (with $\left\Vert \cdot \right\Vert _{\infty }$ referring to the maximum norm) the inequalities
\begin{eqnarray*}
d\left( x,C_{\left( u,b\right) }(t)\right) &\leq &\widetilde{L}\left(
t\right) \left( \left\Vert u(t)-u(s)\right\Vert _{\infty }\left\Vert
x\right\Vert +\left\Vert b(s)-b(t)\right\Vert _{\infty }\right) \\
&\leq &\widetilde{L}\left( t\right) \left( \left\Vert x\right\Vert +1\right)
K\left\vert s-t\right\vert \quad \forall x\in C_{\left( u,b\right) }(s),
\end{eqnarray*}
where $K$ is a Lipschitz constant of $\left( u,b\right)$.
Therefore, if the function $\widetilde{L}\left( t\right) $ is bounded from above on $\left[ 0,T\right]$, say by $L^{\ast }$, then the  desired estimate (\ref{weakhaus0}) would follow with
the function $L\left( \tau \right) :=\left( \tau +1\right) L^{\ast }$, which is 
clearly monotonically increasing. Unfortunately, even for Lipschitzian  controls $\left(
u,b\right)$, the function $\widetilde{L}\left( t\right) $ may be \emph{unbounded from above} as can be seen from the following example.\vspace*{-0.05in}
\begin{example}
\label{counter} {\rm In (\ref{movpoly}) put $m:=2$, $\mathcal{H}:=\mathbb{R}^{2}$,
$T:=1$ and define the smooth (hence Lipschitz continuous) control pair
\begin{equation*}
u_{1}\left( t\right) :=\left( 0,1\right) ;\,b_{1}\left( t\right)
:=1;\,u_{2}\left( t\right) :=\left( t,-1\right) ;\,b_{2}\left( t\right) :=0.
\end{equation*}
For $t\in \left( 0,1\right] $, take $x\left( t\right) :=\left(
t^{-3},1\right) $ and observe that
\begin{equation*}
d\left( x\left( t\right) ,C_{\left( u,b\right) }(t)\right)
=t^{-3}-t^{-1}\;\mbox{ and }\;\max_{i=1,\ldots ,m}\left[ \left\langle u_{i}(t),x\left(
t\right) \right\rangle -b_{i}(t)\right] _{+}=t^{-2}-1.
\end{equation*}
It thus follows from (\ref{hoffman}) that $\widetilde{L}\left( t\right) \geq
t^{-1}$ for all $t\in \left( 0,1\right] $. Therefore, the function $\widetilde{L}\left( t\right) $ is unbounded on $\left[ 0,T\right] $}.
\end{example}\vspace*{-0.2in}
\begin{remark}
\label{specialcases} {\rm There are certain special cases in which Hoffman's error bound leads us to a \emph{bounded} function $\widetilde{L}\left( t\right) $ in \eqref{weakhaus0} on the interval $\left[ 0,T\right]$, even for
non-Lipschitzian controls $\left( u,b\right)$. We mention the following:
\begin{enumerate}
\item In the case of a \emph{moving halfspace} (i.e., $m=1$ and $u(t)\ne 0$ for
all $t\in \left[ 0,1\right] $) with a continuous control $u:\left[ 0,T\right]
\rightarrow \mathcal{H}$ and an arbitrary control $b:\left[ 0,T\right] \rightarrow
\mathbb{R}$, we have that
\begin{equation*}
d\left( x,C_{\left( u,b\right) }(t)\right) =\left\Vert u(t)\right\Vert ^{-1}
\left[ \left\langle u(t),x\right\rangle -b(t)\right] _{+}\leq L^{-1}\left[
\left\langle u(t),x\right\rangle -b(t)\right] _{+}
\end{equation*}
for all $ t\in \left[0,1\right]$ and all $x\in \mathcal{H}$, where $L:=\inf\limits_{t\in \left[ 0,1\right] }\left\Vert u(t)\right\Vert >0$.

\item In the case where variable control functions are situated only on the \emph{right-hand side} of \eqref{movpoly}, i.e, when $u\left( t\right)\equiv
u\ne 0$) while $b:\left[ 0,T\right] \rightarrow \mathbb{R}$ is arbitrary, it follows from \cite[Proposition~4.6]{jourzag} that
\begin{equation*}
d\left( x,C_{\left( u,b\right) }(t)\right) \leq L\max_{i=1,\ldots ,m}\left[
\left\langle u_{i}(t),x\right\rangle -b_{i}(t)\right] _{+}\quad \forall t\in %
\left[ 0,T\right] \,\,\forall x\in \mathcal{H}
\end{equation*}
whenever $C_{\left( u,b\right) }(t)\ne\emp$ for all $t\in \left[ 0,T\right] $.
\end{enumerate}}
\end{remark}\vspace*{-0.05in}

\noindent Example~\ref{counter} illustrates the drastic impact of fully
controlled polyhedral moving sets on Hoffman's error bound starting from
dimension two, even for smooth controls. Fortunately, it turns out that---despite the fact that the approach using Hoffman's error bound sketched
above is not viable for our purposes---we may find an \emph{alternative path} based on (\ref{weakhaus0}), in order to reach the desired goal. To support this
idea, let us revisit Example~\ref{counter} and observe that the sweeping process generated by the Lipschitzian control in this example does admit a unique Lipschitzian solution for an arbitrary initial point $x_{0}\in C_{\left(
u,b\right)}(0)$.\vspace*{-0.05in}
\begin{example}
\label{countercalc} {\rm Consider the control pair $\left( u,b\right)$ defined in Example~\ref{counter} and fix an arbitrary initial point $x_{0}\in C_{\left( u,b\right)}(0)$. We subdivide the initial polyhedron as $C_{\left( u,b\right)
}(0)=\Omega _{1}\cup \Omega _{2}$ with the sets
\begin{equation*}
\Omega _{1}:=\left\{ x\in C_{\left( u,b\right) }(0)\;\big|\;x_{2}<x_{1}\right\}\;\mbox{ and }\;
\Omega _{2}:=\left\{ x\in C_{\left( u,b\right) }(0)\big|\;x_{2}\geq
x_{1}\right\}.
\end{equation*}
If $x_{0}\in \Omega _{2}$, then for an arbitrary time $t\in \left( 0,1\right) $ the boundaries of the two controlled halfspaces have no contact with $x_{0}$.
Consequently, $\dot{x}(t)=0$ for all $t\in \left( 0,1\right) $, and hence $
x\left( t\right) =x_{0}$ for all $t\in \left[ 0,1\right] $. In contrast, for
$x_{0}\in \Omega _{1}$ we get
\begin{align*}
x\left( t\right) &=\left\{
\begin{array}{cc}
x_{0} & t\in \left[ 0,t_{1}\right] \\
y\left( t\right) & t\in \left( t_{1},t_{2}\right) \\
\left( 1/t,1\right) & t\in \left[ t_{2},1\right]
\end{array}
\right., \,\, t_{1}=\frac{x_{0,2}}{x_{0,1}},\,\, t_{2}=\left\{
\begin{tabular}{ll}
$\frac{1}{\sqrt{\left\Vert x_{0}\right\Vert ^{2}-1}}$ & if $\left\Vert
x_{0}\right\Vert \geq \sqrt{2}$ \\
$\infty $ & else
\end{tabular}
\right.,\\
&\quad y_{1}\left( t\right)=\frac{\left\Vert x_{0}\right\Vert }{%
\sqrt{1+t^{2}}},\;\mbox{ and }\;y_{2}\left( t\right)=\frac{\left\Vert
x_{0}\right\Vert }{\sqrt{1+t^{2}}}t.
\end{align*}
Here $t_{1}$ denotes the time when the second halfspace (the moving one)
becomes binding for $x_{0}$ for the first time, i.e., when $tx_{0,1}=x_{0,2}$. This gives us the indicated formula for $t_{1}$. For $t<t_{1}$ both halfspaces are nonbinding for $x_{0}$; so $\dot{x}(t)=0$, and hence $x\left( t\right)
=x_{0}$ for all $t\in \left[ 0,t_{1}\right]$. For $t\geq t_{1}$ the second
halfspace is binding. The first halfspace also becomes binding at a certain time $t_{2}>t_{1}$; so we have $x_{2}\left( t\right) =1$ for all $t\in \lbrack t_{2},1]$. Since the second halfspace keeps binding, it follows that
$tx_{1}\left( t\right) =x_{2}\left( t\right) =1$ from where we conclude that $
x_{1}\left( t\right) =1/t$ during this period of time. It remains to
determine the trajectory $x\left( t\right) $ for $t\in \left(
t_{1},t_{2}\right)$, as well as the switching time $t_{2}$. Since in this
interval only the second halfspace is binding, we derive the following
relations from the sweeping dynamics:
\begin{equation*}
-\dot{x}(t)\in N_{C_{\left( u,b\right) }(t)}\left( x(t)\right) =\mathbb{R}%
_{+}\left( t,-1\right) \quad \forall t\in \left( t_{1},t_{2}\right) .
\end{equation*}
Consequently, there exists a function $\lambda \left( t\right) \leq 0$ such
that
\begin{equation*}
\dot{x}_{1}(t)=t\lambda \left( t\right) ;\quad \dot{x}_{2}(t)=-\lambda
\left( t\right) \quad \forall t\in \left( t_{1},t_{2}\right).
\end{equation*}
On the other hand, with the second halfspace being binding, we also have
that $tx_{1}\left( t\right) =x_{2}\left(t\right)$ for all $t\in\lbrack
t_{1},t_{2})$. This tells us therefore that
\begin{equation*}
\dot{x}_{1}(t)=-t\dot{x}_{2}(t)=-\frac{x_{2}\left( t\right) }{x_{1}\left(
t\right) }\dot{x}_{2}(t)\Longleftrightarrow \dot{x}_{1}(t)x_{1}\left(
t\right) +\dot{x}_{2}(t)x_{2}\left( t\right) =0\quad \forall t\in \left(
t_{1},t_{2}\right) .
\end{equation*}
The solution to the latter differential equation is given by $
x_{1}^{2}\left( t\right) +x_{2}^{2}\left( t\right) =C$, where the constant $
C $ can be identified from the fact that $x\left( t_{1}\right) =x_{0}$, which yields $C=\left\Vert x_{0}\right\Vert ^{2}$. Along with the equality $tx_{1}\left(t\right) =x_{2}\left( t\right) $, we identify the function
$y(t)$ indicated in the formula above. Finally, the switching time $t_{2}$
is determined from the relation $y_{2}\left( t_{2}\right) =1$. Observe that for $\left\Vert x_{0}\right\Vert <
\sqrt{2}$ the first halfspace is never binding in the given time interval $
\left[ 0,1\right]$. It is easy to check that the determined solution $x(t)$ is Lipschitz continuous on the entire interval $[0,1]$, and as such it has to be unique due \cite[Theorem~3]{kunzemonteiro}}.
\end{example}\vspace*{-0.2in}

\subsubsection{\bf Uniform Slater condition and weakened Hausdorff
estimate}

\noindent As shown in our subsequent analysis, the reason why
the announced result---that Lipschitzian controls yield Lipschitzian
solutions of the sweeping process---can be maintained in Example~\ref{counter} despite the fact that an argumentation via Hoffman's error bound
does not apply, consists in the fulfillment of an appropriate \emph{constraint
qualification}. Now we introduce this qualification condition, which plays a crucial role not only in establishing existence and stability results presented in what follows, but also in the two last sections of the paper dealing with the verification of the strong convergence of discrete approximations and the derivation of necessary optimality conditions for sweeping optimal control.\vspace*{0.03in}

Here is this easy formulated and natural qualification condition.\vspace*{-0.05in}

\begin{definition}
We say that the moving polyhedron in \eqref{movpoly} generated by the given control pair $(u,b)$ satisfies the {\sc uniform Slater condition} if
\begin{equation}
\forall t\in \left[ 0,T\right] \,\,\exists x\in \mathcal{H}\;\mbox{ such that}\;\left\langle
u_{i}\left( t\right) ,x\right\rangle <b_{i}\left( t\right) \quad \forall
i=1,\ldots ,m.  \label{unifslater}
\end{equation}
\end{definition}\vspace*{-0.03in}

We emphasize that, unlike the boundedness of $\Tilde L(t)$ in Hoffman's error bound estimate \eqref{hoffman}, this
constraint qualification is \emph{essential} for our desired result.
Indeed, a simple two-dimensional example taken from \cite[Example~2.3]{mor2} shows that, even for smooth control functions, the sweeping process (\ref{sweeping}) may not admit a solution when (\ref{unifslater}) is violated. On the other hand,
we see below that (\ref{unifslater}) yields the weakened Hausdorff
estimate (\ref{weakhaus0}), which is  the first step mentioned in the introduction to this section.\vspace*{0.05in}

Before deriving (\ref{weakhaus0}) via \eqref{unifslater}, we show that the following seemingly stronger version of (\ref{unifslater}) has been used in the earlier work on the existence of solutions to sweeping processes defined by moving polyhedra \cite[Assumption
(H4)]{mor1}:
\begin{equation}
\exists \varepsilon >0\,\,\forall t\in \left[ 0,T\right] \,\,\exists x\in
\mathcal{H}\;\mbox{with}\;\left\langle u_{i}\left( t\right) ,x\right\rangle \leq
b_{i}\left( t\right) -\varepsilon \quad \forall i=1,\ldots ,m
\label{slater2}
\end{equation}
It turns out, however, that this ``strong uniform Slater condition" is \emph{equivalent} to the uniform Slater condition formulated in \eqref{unifslater}.\vspace*{-0.05in}
\begin{proposition}
\label{slaterequiv} Assume that the control $\left( u,b\right)$ in \eqref
{movpoly} is continuous. Then conditions \eqref{unifslater} and \eqref
{slater2} are equivalent.
\end{proposition}\vspace*{-0.15in}
\begin{proof}
Since (\ref{slater2}) obviously yields (\ref{unifslater}), it remains to verify the opposite implication. Assume that  (\ref{slater2}) fails, which tells us that
\begin{equation*}
\forall n\in \mathbb{N}\,\,\exists t_{n}\in \left[ 0,T\right] \,\,\forall
x\in \mathcal{H}\,\,\exists i\in \left\{ 1,\ldots ,m\right\}\;\mbox{with}\;\left\langle
u_{i}\left( t_{n}\right) ,x\right\rangle >b_{i}\left( t_{n}\right) -\frac{1}{
n}.
\end{equation*}
For some subsequence $t_{n_{k}}\in \left[ 0,T\right]$, there exists $\bar{t}%
\in \left[ 0,T\right] $ such that $t_{n_{k}}\rightarrow _{k}\bar{t}$. Fix an
arbitrary vector $x\in \mathcal{H}$ and then get
\begin{equation*}
\forall k\in \mathbb{N}\,\,\exists i_{k}\in \left\{ 1,\ldots ,m\right\}\;\mbox{with}\;\left\langle u_{i_{k}}\left( t_{n_{k}}\right) ,x\right\rangle
>b_{i_{k}}\left( t_{n_{k}}\right) -\frac{1}{n_{k}}.
\end{equation*}
Selecting another subsequence, find $i^{\ast }\in \left\{1,\ldots ,m\right\}$ such that $i_{k_{l}}\equiv i^{\ast }$. Therefore, we have the inequalities
\begin{equation*}
\left\langle u_{i^{\ast }}\left(
t_{n_{k_{l}}}\right) ,x\right\rangle >b_{i^{\ast }}\left(
t_{n_{k_{l}}}\right) -\frac{1}{n_{k_{l}}}\;\mbox{ for all }\;l\in \mathbb{N}.
\end{equation*}
Passing  there to the limit as $l\rightarrow \infty $ gives us $\left\langle
u_{i^{\ast }}\left( \bar{t}\right) ,x\right\rangle \geq b_{i^{\ast }}\left(
\bar{t}\right) $. Since $x\in \mathcal{H}$ was chosen arbitrarily, we arrive at 
\begin{equation*}
\exists \bar{t}\in \left[ 0,T\right] \,\,\forall x\in \mathcal{H}\,\,\exists
i^{\ast }\in \left\{ 1,\ldots ,m\right\}\;\mbox{with}\;\left\langle u_{i^{\ast }}\left(
\bar{t}\right) ,x\right\rangle \geq b_{i^{\ast }}\left( \bar{t}\right),
\end{equation*}
which contradicts (\ref{unifslater}) and thus completes the proof of the proposition.\qed
\end{proof}

\noindent Now we turn to the announced proof of the weakened Hausdorff
estimate (\ref{weakhaus0}). Given $\dd>0$, define
the \emph{$\delta -$moving polyhedron} by
\begin{equation}
C_{(u,b)}^{\left( \delta \right) }(t):=\big\{ x\in \mathcal{H}\;\big|\;\left\langle
u_{i}(t),x\right\rangle \leq b_{i}(t)-\delta \,\, \left(i=1,\ldots
,m\right)\big\} \,\,\left( t\in \left[ 0,T\right] \right) .
\label{delmov}
\end{equation}

To proceed, we first present the following crucial technical lemma involving continuous controls $(u,b)\in\mathcal{C}([0,T],\mathcal{H}^{m})\times \mathcal{C}([0,T],\mathbb{
R}^{m})$ in the moving polyhedron \eqref{movpoly} endowed with the maximum norm
\begin{equation*}
\left\Vert (u,b)\right\Vert _{\infty }:=\max_{t\in \left[ 0,T\right]
,i=1,\ldots ,m}\left\Vert u_{i}(t)\right\Vert +\max_{t\in \left[ 0,T\right]
,i=1,\ldots ,m}\left\vert b_{i}(t)\right\vert.
\end{equation*}
The associated closed ball in this space centered at $(u,b)$ with radius $r>0$ is denoted by $\mathbb{B}_{\infty }\left( (u,b),r\right)$.\vspace*{-0.05in}
\begin{lemma}
\label{strongselection} Fix continuous control $(\bar{u},
\bar{b})\in \mathcal{C}([0,T],\mathcal{H}^{m})\times \mathcal{C}([0,T],
\mathbb{R}^{m})$ satisfying the uniform Slater condition \eqref{unifslater}. Then there exists $\varepsilon>0$ such that whenever $\gamma\in \left( 0,\varepsilon\right)$ we can find a continuous function $\widehat{x}\in \mathcal{C}([0,T],\mathcal{H})$ for which
\begin{equation}
\widehat{x}(t)\in C_{(u,b)}^{\left( \gamma \right) }(t)\,\, \forall t\in %
\left[ 0,T\right] \,\,\forall (u,b)\in \mathcal{B}:=\mathbb{B}_{\infty
}\left( \left( \bar{u},\bar{b}\right),\frac{\varepsilon -\gamma}{3\left(
1+\left\Vert \widehat{x}\right\Vert _{\infty }\right) }\right).
\label{select1}
\end{equation}
Furthermore, we have the estimate
\begin{equation}
d(x,C_{(u,b)}(t))\leq {\frac{{f_{(u,b)}(t,x)}}{{f_{(u,b)}(t,x)-f_{(u,b)}(t,
\widehat{x}(t))}}}\Vert x-\widehat{x}(t))\Vert\,\,\forall t\in \left[ 0,T\right] \label{select2}
\end{equation}
for all $t\in \left[ 0,T\right]$, all $x\in \mathcal{H}\backslash C_{(u,b)}(t)$, and all $(u,b)\in \mathcal{B}$, where $f_{(u,b)}(t,x):=\max_{i=1,\cdots ,m}\langle u_{i}(t),x\rangle
-b_{i}(t)$. Finally,
\begin{align}
&d(x,C_{(u^{\prime },b^{\prime })}(t))\leq\notag \\
&\Vert x-\widehat{x}(t)\Vert \min
\left\{ 1,\gamma ^{-1}\max_{i=1,\cdots ,m}\left[ \langle u_{i}^{\prime
}(t)-u_{i}(s),x\rangle +b_{i}(s)-b_{i}^{\prime }(t)\right] _{+}\right\}
\label{select3}
\end{align}
for all $(u,b),(u^{\prime },b^{\prime })\in \mathcal{B}$, all $s,t\in
\left[ 0,T\right], \,$\ and all $x\in C_{(u,b)}(s)$.
\end{lemma}\vspace*{-0.15in}
\begin{proof} As shown in Proposition~\ref{slaterequiv}, the imposed uniform Slater condition (\ref{unifslater}) is equivalent to 
(\ref{slater2}) for $(u,b):=(\bar{u},\bar{b})$. Using the latter and choosing $\varepsilon >0$ therein, pick an arbitrary number $\gamma \in \left( 0,\varepsilon \right)$ and define
\begin{equation*}
\delta :=\frac{2\varepsilon +\gamma }{3}\in \left( 0,\varepsilon \right) .
\end{equation*}
Then condition (\ref{slater2}) tells us that
\begin{equation*}
\forall t\in \left[ 0,T\right] \,\,\exists x\in \mathcal{H}\;\mbox{ with }\;\left\langle
\bar{u}_{i}\left( t\right) ,x\right\rangle \leq \bar{b}_{i}\left( t\right)
-\varepsilon <\bar{b}_{i}\left( t\right) -\delta \quad \forall i=1,\ldots,m.
\end{equation*}
In other words, for each $t\in \left[ 0,T\right]$ the convex set $C_{\left(
\bar{u},\bar{b}\right) }^{\left( \delta \right) }(t)$ admits a Slater point.
This ensures the inclusion
\begin{equation*}
C_{\left( \bar{u},\bar{b}\right) }^{\left( \delta \right) }(t)\subseteq
\mathrm{cl}\,\left\{ x\in \mathcal{H}\;\big|\;\left\langle \bar{u}_{i}\left(
t\right) ,x\right\rangle <\bar{b}_{i}\left( t\right) -\delta \right\} \quad
\forall t\in \left[ 0,T\right]
\end{equation*}
which in turn allows to conclude (by invoking, e.g., 
\cite[Theorem~3.1.5]{fivemen}) that $C_{( u,b)}^{(\delta)}:[0,T]\rightrightarrows \mathcal{H}$ is a lower
semicontinuous multifunction. Since the images $C_{(\ou,\ob)}^{(\dd)}(t)$ are closed and convex for all $t\in \left[0,T\right]$, the classical Michael selection theorem ensures the existence of a continuous function $\widehat{x}\in \mathcal{C}([0,T],\mathcal{H})$
with
\begin{equation*}
\widehat{x}(t)\in C_{\left( \bar{u},\bar{b}\right) }^{\left( \delta \right)
}\left( t\right) \quad \forall t\in \left[ 0,T\right].
\end{equation*}
Next we fix an arbitrary continuous control $(u,b)\in \mathcal{B}$ and get by the definition of $\delta $ the following inequalities:
\begin{eqnarray*}
\left\langle u_{i}\left( t\right) ,\widehat{x}(t)\right\rangle -b_{i}\left(
t\right) &\leq &\left\langle \bar{u}_{i}\left( t\right) ,\widehat{x}
(t)\right\rangle +\left\Vert u_{i}\left( t\right) -\bar{u}_{i}\left(
t\right) \right\Vert\cdot\left\Vert \widehat{x}(t)\right\Vert -b_{i}\left(
t\right) \\
&\leq &\bar{b}_{i}\left( t\right) -\delta +\left\Vert u_{i}\left( t\right) -%
\bar{u}_{i}\left( t\right) \right\Vert\cdot\left\Vert \widehat{x}(t)\right\Vert
-b_{i}\left( t\right) \\
&\leq &\frac{2}{3}\left( \varepsilon -\gamma \right) -\delta \leq -\gamma
\quad \forall t\in \left[ 0,T\right] \,\,\forall i=1,\ldots ,m.
\end{eqnarray*}
Thus $\widehat{x}\in \mathcal{C}([0,T],\mathcal{H})$ and $\widehat{x}(t)\in
C_{\left( u,b\right) }^{\left( \gamma \right) }\left( t\right) $ for all $
t\in \left[ 0,T\right]$, which verify (\ref{select1}). 

Addressing the second assertion of the lemma, fix arbitrary elements $t\in \left[ 0,T\right] $, $(u,b)\in
\mathcal{B}$, and $x\in \mathcal{H}\backslash C_{(u,b)}(t)$. Remembering the construction of ${f_{(u,b)}}$, we have
that ${f_{(u,b)}(t,x)>0}$ by $x\in \mathcal{H}\backslash C_{(u,b)}(t)$ and ${
f_{(u,b)}(t,\widehat{x}(t))\leq -\gamma <0}$ by the already proved relation (
\ref{select1}), define
\begin{equation*}
\lambda :={\frac{{f_{(u,b)}(t,x)}}{{f_{(u,b)}(t,x)-f_{(u,b)}(t,\widehat{x}
(t))}}\in }\left( 0,1\right).
\end{equation*}
It follows from the convexity of ${f_{(u,b)}}(t,\cdot )$ that
\begin{equation*}
f_{(u,b)}(t,(1-\lambda )x+\lambda \widehat{x}(t))\leq (1-\lambda
)f_{(u,b)}(t,x)+\lambda f_{(u,b)}(t,\widehat{x}(t))=0,
\end{equation*}
and so $(1-\lambda )x+\lambda \widehat{x}(t)\in C_{(u,b)}(t)$. This verifies (\ref{select2}), which can be written as
\begin{equation*}
d(x,C_{(u,b)}(t))\leq \Vert x-((1-\lambda )x+\lambda \widehat{x}(t))\Vert ={
\lambda }\Vert x-\widehat{x}(t)\Vert.
\end{equation*}

It remains to justify the final assertion of the lemma. To proceed, fix arbitrary elements $s,t\in \lbrack 0,T]$, $
(u,b),(u^{\prime },b^{\prime })\in \mathcal{B}$, and $x\in C_{(u,b)}(s)$. If $
x\in C_{(u^{\prime },b^{\prime })}(t)$, then (\ref{select3}) holds
trivially. Supposing now that $x\notin C_{(u^{\prime },b^{\prime })}(t)$ gives us ${
f_{(u^{\prime },b^{\prime })}(t,x)>0}$ and ${f_{(u,b)}(t,\widehat{x}(t))\le
-\gamma }$ by (\ref{select1}). Therefore, (\ref{select2}) yields
\begin{align*}
&d(x,C_{(u^{\prime },b^{\prime })}(t))\\&\leq{\frac{{f_{(u^{\prime
},b^{\prime })}(t,x)}}{{f_{(u^{\prime },b^{\prime })}(t,x)-f_{(u^{\prime
},b^{\prime })}(t,\widehat{x}(t))}}}\Vert x-\widehat{x}(t))\Vert
\leq {\gamma }^{-1}{f_{(u^{\prime },b^{\prime })}(t,x)}\Vert x-\widehat{x}(t))\Vert
\\
& \leq {\gamma }^{-1}\left( {f_{(u^{\prime },b^{\prime
})}(t,x)-f_{(u,b)}(s,x)}\right) \Vert x-\widehat{x}(t))\Vert \quad \left(
\mbox{because of }x\in C_{(u,b)}(s)\right) \\
& \leq {\gamma }^{-1}\Vert x-\widehat{x}(t)\Vert \max_{i=1,\cdots ,m}\left[
\langle u_{i}^{\prime }(t)-u_{i}(s),x\rangle +b_{i}(s)-b_{i}^{\prime }(t)
\right] _{+}.
\end{align*}
Since $\widehat{x}(t)\in C_{(u^{\prime },b^{\prime })}^{\left( \gamma
\right) }(t)\subseteq C_{(u^{\prime },b^{\prime })}(t)$ by (\ref{select1}),
we also have that $d(x,C_{(u^{\prime },b^{\prime })}(t))\leq \Vert x-
\widehat{x}(t))\Vert $. Combining the above verifies (\ref{select3}) and  completes the proof. \qed
\end{proof}

We are now in a position to derive the weakened Hausdorff estimate (\ref{weakhaus0}).\vspace*{-0.05in}

\begin{theorem}\label{slaterest} Let $\left( u,b\right) $ be a Lipschitz continuous control along which the moving polyhedron \eqref{movpoly} satisfies the
uniform Slater condition \eqref{unifslater}. Then there exist constants $
K_{1},K_{2}\geq 0$ such that the weakened Hausdorff estimate \eqref{weakhaus0} holds with the  monotonically increasing function $L:\mathbb{R}_{+}\rightarrow \mathbb{R}_{+}$  defined by
\begin{equation}
L\left( r\right):=K_{1}\left( r+1\right) \left( r+K_{2}\right) \quad \left(
r\geq 0\right).  \label{lrquad}
\end{equation}
\end{theorem}\vspace*{-0.15in}
\begin{proof} We again employ the uniform Slater condition \eqref{unifslater} in the equivalent form \eqref{slater2} by Proposition~\ref{slaterequiv}. Then we get from \eqref{select3} in Lemma~\ref{strongselection} that
\begin{equation*}
d(x,C_{(u,b)}(t))\leq \frac{2}{\varepsilon }\Vert x-\widehat{x}(t)\Vert
\max_{i=1,\cdots ,m}\left[\langle u_{i}(t)-u_{i}(s),x\rangle
+b_{i}(s)-b_{i}(t)\right] _{+}
\end{equation*}
along a continuous function $\widehat{x}(\cdot)$ for all $s,t\in\left[ 0,T
\right] \,$\ and all $x\in C_{(u,b)}(s)$. Define $\varkappa
:=\max\limits_{t\in \left[ 0,T\right] }\left\Vert \widehat{x}\left( t\right)
\right\Vert \geq 0$ and denote by $K\geq 0$ a Lipschitz constant of the
control pair $\left( u,b\right) $. Then we have the estimate
\begin{equation*}
d(x,C_{(u,b)}(t))\leq \frac{2K}{\varepsilon }\Vert x-\widehat{x}(t)\Vert
\left( \left\Vert x\right\Vert +1\right) \left\vert s-t\right\vert \leq
\frac{2K}{\varepsilon }\left( \Vert x\Vert +\varkappa \right) \left(
\left\Vert x\right\Vert +1\right) \left\vert s-t\right\vert
\end{equation*}
for all $s,t\in \left[ 0,T\right] \,$\ and all $x\in C_{(u,b)}(s)$. This is
exactly (\ref{weakhaus0}) with the monotonically increasing function $
L\left( r\right) :=\delta ^{-1}K\left( r+\varkappa \right) \left( r+1\right)
$. \qed
\end{proof}\vspace*{-0.2in}
\begin{remark}
\noindent {\rm The moving polyhedron $C_{\left( u,b\right) }$ defined in Example~\ref{counter} does satisfy the uniform Slater condition. To see this, select the constant solution $x\left( t\right) \equiv \left(0,0.5\right)$
in (\ref{unifslater}). Thus the estimate (\ref{weakhaus0}) can be verified
in this example via Theorem~\ref{slaterest}, while the usage of
Hoffman's error bound does not lead us to the desired result. The reason is
that Hoffman's error bound---if applicable as in the special cases mentioned
in Remark~\ref{specialcases}---would necessarily bring us to an \emph{affine function} $L$ in \eqref{weakhaus0}; see the discussion above in 
Example~\ref{counter}. Yet, a closer inspection of the example shows that such an affine function $L$ cannot work in this example. Indeed, consider the
sequences
\begin{equation*}
x^{\left( n\right) }:=\left( 2n,0\right) \in C_{\left( u,b\right) }\left(
0\right) ;\quad t_{n}:=n^{-1}\quad \left( n\in \mathbb{N}\right).
\end{equation*}
Assuming that estimate (\ref{weakhaus0}) holds with an affine function $
L\left( r\right):=ar+b$ and choosing $s:=0$, we arrive at the following contradiction 
\begin{eqnarray*}
n &\leq &\sqrt{1+n^{2}}=d\left( x^{\left( n\right) },C_{\left( u,b\right)
}(t_{n})\right) \leq \left( a\left\Vert x^{\left( n\right) }\right\Vert
+b\right) t_{n} \\
&=&\left(2an+b\right) n^{-1}\leq 2a+\left\vert b\right\vert \quad \forall
n\in \mathbb{N}.
\end{eqnarray*}
On the other hand, the choice of the \emph{quadratic
function} (\ref{lrquad}) by Theorem~\ref{slaterest} allows us to derive the weakened Hausdorff estimate \eqref{weakhaus0} in this example}.
\end{remark}\vspace*{-0.2in}

\subsubsection{\bf General truncation lemma}

\noindent The last subsection of this section accomplishes the \emph{second step} of our approach outlined in the introduction to this section. The following general truncation result clearly implies the desired estimates \eqref{truncest} for truncating polyhedra.\vspace*{-0.05in} 
\begin{lemma}
\label{TL} Let $(X,\Vert \cdot \Vert )$ be a normed space, and let $C$ 
be a nonempty, closed, and convex subset of $X$. Define the truncating 
set $C^{r}:=C\cap \mathbb{B}\left(0,r\right) $ for $r>0$. Then we have the estimate
\begin{equation}
d(x,C^{r})\leq \frac{2r}{r-d(0, C)} d(x,C)\quad \forall x\in \mathbb{B}
\left( 0,r\right) \,\,\forall r>d\left( 0,C\right).  \label{truncest0}
\end{equation}
Consequently, it follows that
\begin{equation}
d(x,C^{r})\leq 3d(x,C)\quad \forall x\in \mathbb{B}\left( 0,r\right)
\,\,\forall r> 3d\left( 0,C\right).  \label{truncest1}
\end{equation}
\end{lemma}\vspace*{-0.15in}
\begin{proof}
Pick arbitrary elements $r>d\left( 0,C\right) $, $x\in \mathbb{B}\left( 0,r\right) $, and $\varepsilon $ with $0<\varepsilon <r-d(0,C)$. If $x\in C$, then
$x\in C^{r}$ and (\ref{truncest0}) holds trivially. Assume now
that $x\notin C$, and so $d(x,C)>0$. Choose $x_{0},y\in C$ such that
\begin{equation}
\left\Vert x_{0}\right\Vert \leq \beta :=d\left( 0,C\right) +\varepsilon
,\,\,\left\Vert x-y\right\Vert \leq d(x,C)+\min \left\{ \varepsilon ,d\left(
x,C\right) \right\} .  \label{tworel}
\end{equation}
If $\Vert y\Vert \leq r$, then $y\in C^{r}$, and (\ref{truncest0}) follows
from the inequality in (\ref{tworel}). Therefore, it remains to examine the case where $\Vert y\Vert >r$. The equality in (\ref{tworel}) combined with $
\varepsilon <r-d(0,C)$ gives us the estimate $\left\Vert x_{0}\right\Vert \leq \beta <r$. Therefore, there exists $\gamma \in \left( 0,1\right)$ such that $\Vert
z\Vert =r$ for $z:=(1-\gamma )y+\gamma x_{0}$. The convexity of $C$ readily ensures that $z\in C^{r}$. Then we have
\begin{equation*}
r\leq (1-\gamma )\Vert y\Vert +\gamma \Vert x_{0}\Vert \text{\quad
\mbox{or,
equivalently,}\quad }\gamma \left( \Vert y\Vert -\Vert x_{0}\Vert \right)
\leq \Vert y\Vert -r.
\end{equation*}
Due to $\Vert y\Vert >r>\beta \geq \left\Vert x_{0}\right\Vert $, the latter implies that
\begin{equation*}
\Vert z-y\Vert =\gamma \Vert y-x_{0}\Vert \leq \frac{\Vert y\Vert -r}{\Vert
y\Vert -\beta }\left( \Vert y\Vert +\beta \right).
\end{equation*}
Taking into account that $\Vert x\Vert \leq r$ brings us to
\begin{equation*}
\Vert y\Vert \leq \Vert y-x\Vert +\Vert x\Vert \leq d(x,C)+\varepsilon +r,
\end{equation*}
and therefore we arrive at the estimate
\begin{equation*}
\Vert z-y\Vert \leq \frac{\Vert y\Vert +\beta }{\Vert y\Vert -\beta }
(d(x,C)+\varepsilon).
\end{equation*}
Combining all the above leads us to the relationships
\begin{equation*}
\Vert z-x\Vert \leq \Vert z-y\Vert +\Vert y-x\Vert \leq (1+\frac{\Vert
y\Vert +\beta }{\Vert y\Vert -\beta })(d(x,C)+\varepsilon )\leq \left( 2+%
\frac{2\beta }{r-\beta }\right) (d(x,C)+\varepsilon ).
\end{equation*}%
Since $z\in C^{r}$ and $\varepsilon$ was chosen arbitrarily with $
0<\varepsilon <r-d(0,C)$, we get
\begin{equation*}
d(x,C_{r})\leq \left( 2+\frac{2d(0, C) }{r-d(0, C) }\right) d(x,C),
\end{equation*}
which verifies (\ref{truncest0}) and thus completes the proof of the truncation lemma.  \qed
\end{proof}\vspace*{-0.25in}

\section{Existence and uniqueness of sweeping solutions}\label{exsol}\vspace*{-0.05in}

\noindent The main goal of this section is establishing two class-preservation \emph{existence and uniqueness} theorems for polyhedral controlled sweeping processes defined in \eqref{movpoly} and \eqref{sweeping} under the uniform Slater condition \eqref{unifslater} in the setting of \emph{separable Hilbert spaces}. Namely, we aim at proving that \emph{Lipschitz continuous} controls $(u,b)$ uniquely generate Lipschitz continuous trajectories of ${\cal S}_{(u,b)}$ and that \emph{absolutely continuous} (of class $W^{1,1}$) controls uniquely generate sweeping trajectories of the same class. Note that results of this type in the $W^{1,2}$ control-trajectory framework we obtained in  \cite{tolsto,tolsto2,tolsto3} for various types of sweeping processes under appropriate assumptions in separable Hilbert spaces. Similar preservation results of class $W^{1,2}$ were established in \cite{mor1} in finite dimensions under the strong uniform Slater condition \eqref{slater2} reducing to \eqref{unifslater} as we now know. Observe also that results of this type in class of $W^{1,1}$ were derived in \cite{mor2,colmor} for polyhedral sweeping processes in finite-dimensional spaces under essentially stronger qualification conditions than \eqref{unifslater} used in what follows. Our approach below is strongly based on the truncation procedure and error bound estimates developed in the previous section.\vspace*{0.03in} 

Here is the first theorem dealing with Lipschitzian controls.\vspace*{-0.08in} 

\begin{theorem}
\label{existlip} Let $\mathcal{H}$ be a separable Hilbert space. Assume that $\left( u,b\right)$ is Lipschitz continuous control and
that the moving polyhedron $C_{\left(u,b\right)}$ in \eqref{movpoly}
satisfies the uniform Slater condition \eqref{unifslater} along this control pair. Then the sweeping process $\left( \mathcal{S}_{\left( u,b\right) }\right)$ admits a unique Lipschitz continuous solution.
\end{theorem}\vspace*{-0.2in}
\begin{proof} Theorem~\ref{slaterest} ensures the existence of a monotonically
increasing function $L:\mathbb{R}_{+}\rightarrow \mathbb{R}_{+}$ satisfying the weakened Hausdorff estimate \eqref{weakhaus0}. This gives us for each $r>0$ a constant $\widehat{L}_{r}:=L\left( r\right)$ such that (\ref{weakhaus}) holds. Thus for all $r>0$, all $s,t\in \left[ 0,T\right] $, and all $x\in C_{\left( u,b\right) }(s)$ with $\left\Vert x\right\Vert\leq r$ there is
$y\in C_{\left( u,b\right) }(t)$ satisfying
\begin{equation*}
\left\Vert x-y\right\Vert \leq \left( \widehat{L}_{r}+1\right) \left\vert
s-t\right\vert.
\end{equation*}
Indeed, the latter is obvious with the choice of $y:=x$ in the case where $s=t$, and this follows from (\ref{weakhaus}) and from $d\left( x,C_{\left( u,b\right)
}(t)\right) <\left( \widehat{L}_{r}+1\right) \left\vert s-t\right\vert $ in
the case where $s\neq t$. Since the linear function $s\longmapsto \left( 
\widehat{L}_{r}+1\right) s$ trivially belongs to $W^{1,2}\left[ 0,T\right]$, it is $r$-\textit{weakly uniformly lower semicontinuous from the right}
for $p=2$ in the sense of Tolstonogov \cite[eq. (2.2)]{tolsto}. Therefore, we deduce from \cite[Lemma~2.1 and Lemma~3.1]{tolsto} that the sweeping
process $\left( \mathcal{S}_{\left( u,b\right) }\right)$ has a unique
solution $x^{\ast }\in W^{1,2}\left( \left[ 0,T\right],\mathcal{H}\right)$.
In particular, the trajectory $x^{\ast}(t)$ is absolutely continuous on $[0,T]$. It remains to show that $x^{\ast}(t)$ is Lipschitz continuous on this interval. To proceed, define 
\begin{equation}
\rho :=\max\limits_{t\in \left[ 0,T\right] }\left\Vert x^{\ast }\left(
t\right) \right\Vert ;\quad r:=3\rho +1  \label{rhordef}
\end{equation}
and then fix arbitrary $s,t\in \left[ 0,T\right] $ and 
\begin{equation*}
x\in C_{\left( u,b\right)}^{r}(s):=C_{\left( u,b\right) }(s)\cap \mathbb{B}
\left( 0,r\right).
\end{equation*}
As a solution to $\left( \mathcal{S}_{\left( u,b\right) }\right) $, the function $x^{\ast
}(t)$ satisfies the hidden state constraint $x^{\ast }\left( t\right)
\in C_{\left( u,b\right) }(t)$. Therefore, we obtain
\begin{equation*}
r=3\rho +1\geq 3\left\Vert x^{\ast }\left( t\right) \right\Vert +1>3d\left(
0,C_{\left( u,b\right) }(t)\right).
\end{equation*}
This allows us to invoke the truncation result from Lemma~\ref{TL} to get
\begin{equation}
d\left( x,C_{\left( u,b\right) }^{r}(t)\right) \leq 3d\left( x,C_{\left(
u,b\right) }(t)\right).\label{almost}
\end{equation}
On the other hand, Theorem~\ref{slaterest} yields (\ref{weakhaus0}) and
hence gives us a constant $\widehat{L}$ such that (\ref{weakhaus}) holds for our selected $s,t\in \left[ 0,T\right] $. Combining this with (\ref{almost}),
and recalling that $s,t,x$ were chosen arbitrarily, we arrive at the estimate
\begin{equation*}
d\left( x,C_{\left( u,b\right) }^{r}(t)\right) \leq 3\widehat{L}\left\vert
s-t\right\vert \quad \forall s,t\in \left[ 0,T\right] \,\,\forall x\in
C_{\left( u,b\right) }^{r}(s).
\end{equation*}
Interchanging the roles of $s$ and $t$ readily yields the desired Lipschitz
Hausdorff estimate (\ref{haustrunc}) of the truncated moving polyhedron
with modulus $3\widehat{L}$. Employing the standard existence result from \cite[Theorem~2]{kunzemonteiro}) leads us to deducing from the obtained estimate that the truncated sweeping process $\big(\mathcal{\widetilde{S}}_{\left(u,b\right)}\big)$ defined as
\begin{equation}
-\dot{x}
(t)\in N_{C_{\left( u,b\right)}^{r}(t)}\left( x(t)\right)\;\mathrm{
a.e.\,}\;t\in[0,T],\;\; x(0)=x_{0}\in C_{\left( u,b\right)
}^{r}(0) \label{tildes}
\end{equation}
admits a Lipschitz continuous solution $\widetilde{x}(\cdot)$. It follows from the definitions in (\ref{rhordef}) that for all $r>\rho $ we have the inclusions
\begin{equation*}
x^{\ast }\left( t\right) \in C_{\left( u,b\right) }(t)\cap \mathbb{B}\left(
0,\rho \right) \subseteq C_{\left( u,b\right) }(t)\cap \mathrm{int\,}\mathbb{
B}\left( 0,r\right) \subset C_{\left( u,b\right) }^{r}(t)\quad \forall t\in 
\left[ 0,T\right].
\end{equation*}
On the one hand, the resulting inclusion justifies the feasibility of the initial point in $\big(\mathcal{\widetilde{S}}_{\left( u,b\right) }\big)$
due to $x_{0}=x^{\ast }\left( 0\right)$. On the other hand, it tells us that
\begin{equation*}
N_{C_{\left( u,b\right) }^{r}(t)}\left(x^{\ast }(t)\right) =N_{C_{\left(
u,b\right) }(t)}\left(x^{\ast }(t)\right) \quad \forall \mathrm{\,}t\in 
\left[ 0,T\right].
\end{equation*}
Therefore, $x^{\ast }(\cdot)$ being a solution to $\big(\mathcal{S}_{\left(
u,b\right) }\big)$ is also a solution to $\big( \mathcal{\widetilde{S}}
_{\left( u,b\right) }\big)$. Since $x^{\ast}(t)$ is absolutely
continuous on $[0,T]$ as an element of $W^{1,2}\left( \left[ 0,T\right] ,\mathcal{H}\right)$, and since $\left( \mathcal{\widetilde{S}}_{\left( u,b\right) }\right) $ can have at most one absolutely continuous solution by \cite[Theorem~3]{kunzemonteiro}, we conclude that $x^{\ast}(\cdot)=\widetilde{x}(\cdot)$. This ensures that $x^{\ast
}(t)$ is Lipschitz continuous on $[0,T]$, since $\widetilde{x}(t)$ is so. Thus we complete the proof.\qed
\end{proof}\vspace*{-0.05in}

Our next goal in this section is establish the existence of a unique \emph{absolutely continuous} solution of the sweeping process $\big(\mathcal{S}_{\left(u,b\right)}\big)$ generated by any absolutely control $(u,b)$ in the moving polyhedron \eqref{movpoly} under the same uniform Slater condition. Recall that the norms on the spaces of
absolutely continuous functions $W^{1,1}([0,T],\mathcal{H}^{m})$ and $
W^{1,1}([0,T],\mathbb{R}^{m})$ are defined, respectively, by 
\begin{equation*}
\Vert u\Vert _{1,1}:=\sum_{i=1}^{m}\Vert u_{i}(0)\Vert
+\sum_{i=1}^{m}\int_{0}^{T}\Vert \dot{u}_{i}(t)\Vert dt,\;\Vert
b\Vert _{1,1}:=\sum_{i=1}^{m}|b_{i}(0)|+\sum_{i=1}^{m}\int_{0}^{T}|\dot{b}
_{i}(t)|dt.
\end{equation*}
The norm on the product space $W^{1,1}([0,T],\mathcal{H}^{m})\times
W^{1,1}([0,T],\mathbb{R}^{m})$) is $\Vert (u,b)\Vert
_{1,1}:=\Vert u\Vert _{1,1}+\Vert b\Vert _{1,1}$, and the induced ball around $
\left( u,b\right) $ with radius $r$ is $\mathbb{B}_{1,1}\left( \left(
u,b\right),r\right)$.\vspace*{0.05in}

The proof of the following theorem elaborates a reduction idea from \cite{thibault} that allows us to deal with non-Lipschitzian controls of the sweeping dynamics.\vspace*{-0.05in}
\begin{theorem}
\label{existsobolev} Let $\mathcal{H}$ be a separable Hilbert space. Take $\left( \bar{u},\bar{b}\right)\in W^{1,1}([0,T],\mathcal{
H}^{m})\times W^{1,1}([0,T],\mathbb{R}^{m})$ and suppose that the moving
polyhedron $C_{\left( u,b\right) }$ in \eqref{movpoly} satisfies the uniform
Slater condition \eqref{unifslater}. Then the control pair $(u,b)$ generates a unique solution $x\in W^{1,1}\left( \left[ 0,T\right] ,\mathcal{H}\right)$ of the sweeping process $\big(\mathcal{S}_{\left(
u,b\right)}\big) $ in \eqref{sweeping}.
\end{theorem}\vspace*{-0.2in}
\begin{proof}
It follows from the Newton-Leibniz formula that
\begin{equation*}
\left\Vert f\left( t\right) -f\left( s\right) \right\Vert \leq
\int_{s}^{t}\left\Vert \dot{f}\left( r\right) \right\Vert dr\quad \forall
f\in W^{1,1}\left( \left[ 0,T\right] ,\mathcal{H}\right) 
\end{equation*}
whenever $s,t\in\left[ 0,T\right]$ with $s\leq t$. Therefore, for all such $s,t$ we have 
\begin{align}
\sum_{i=1}^{m}\Vert u_{i}(t)-u_{i}(s)\Vert +|b_{i}(t)-b_{i}(s)|&\leq
\left\vert \int_{s}^{t}\sum_{i=1}^{m}\Vert \dot{u}_{i}(r)\Vert +|\dot{b}
_{i}(r)|dr+t-s\right\vert\notag\\&=|\gamma (t)-\gamma (s)|  \label{trafo}
\end{align}
with the strongly increasing and absolutely continuous function 
\begin{equation}\label{gamma}
\gamma (t):=t+\int_{0}^{t}\sum_{i=1}^{m}\Vert \dot{u}_{i}(r)\Vert +|\dot{b}
_{i}(r)|dr
\end{equation}
For each index $i=1,\ldots ,m$, introduce the pair $\left( u_{i}^{\prime},b_{i}^{\prime }\right) :[0,\gamma (T)]\rightarrow H\times \mathbb{R}$ by 
\begin{equation*}
\left( u_{i}^{\prime },b_{i}^{\prime }\right) (\tau ):=\left(
u_{i},b_{i}\right) (\gamma ^{-1}(\tau )),\quad \tau \in \lbrack 0,\gamma (T)].
\end{equation*}
Then we readily have the relationship 
\begin{equation}
C_{\left( u^{\prime },b^{\prime }\right) }(\tau )=C_{\left( u,b\right)
}(\gamma ^{-1}(\tau )),\quad \tau \in \lbrack 0,\gamma (T)].  \label{ctrafo}
\end{equation}
Since $\gamma ^{-1}(0)=0$, it follows from \eqref{ctrafo} that $x_{0}\in
C_{\left( u,b\right) }(0)=C_{\left( u^{\prime },b^{\prime }\right) }(0)$. Therefore, the sweeping process
\begin{equation*}
\big( \mathcal{S}_{\left( u^{\prime },b^{\prime }\right) }^{\prime }\big)
:\qquad -\dot{x}(\tau )\in N_{C_{\left( u^{\prime },b^{\prime }\right)
}(\tau )}\left( x(\tau )\right) \quad \mathrm{a.e.\,}\tau \in \left[
0,\gamma (T)\right] ,\quad x(0)=x_{0}
\end{equation*}
is exactly of type $\left( \mathcal{S}_{\left( u,b\right) }\right) 
$ as in (\ref{sweeping}). Furthermore, (\ref{trafo}) yields 
\begin{equation*}
\Vert u_{i}^{\prime }(\tau _{1})-u_{i}^{\prime }(\tau _{2})\Vert
+|b_{i}^{\prime }(\tau _{1})-b_{i}^{\prime }(\tau _{2})|\leq \left\vert \tau
_{1}-\tau _{2}\right\vert \quad \forall \tau _{1},\tau _{2}\in \left[
0,\gamma (T)\right] \,\,\forall i=1,\ldots ,m,
\end{equation*}
which tells us that the control $\left( u^{\prime },b^{\prime }\right) $ is Lipschitz continuous on the interval $[0,\gg(T)]$. Observe also that $C_{\left( u^{\prime },b^{\prime }\right) }$ satisfies
the uniform Slater condition (\ref{unifslater}) on this interval since 
$C_{\left( u,b\right) }$ does so on the original
interval $\left[ 0,T\right] $). This allows us to invoke
Theorem~\ref{existlip}, applied now to the control $\left( u^{\prime },b^{\prime
}\right) $, and conclude that the modified sweeping process $\big( \mathcal{S}_{\left( u^{\prime },b^{\prime
}\right) }^{\prime }\big) $ admits a unique Lipschitzian solution $y(\cdot)$ with some modulus $K$. For all $t\in \lbrack 0,T]$, set $z\left( t\right)
:=y\left( \gamma (t)\right) $, which implies that $\dot{z}\left( t\right) :=\dot{y}
\left( \gamma (t)\right) \dot{\gamma}(t)$ for a.e.\ $t\in[0,T]$. Hence
\begin{equation*}
\left\Vert \dot{z}\left( t\right) \right\Vert \leq \left\Vert \dot{y}\left(
\gamma (t)\right) \right\Vert \dot{\gamma}(t)\leq K\dot{\gamma}(t)\quad 
\mathrm{a.e.\;}\;t\in \left[ 0,T\right] .
\end{equation*}
Since $y(\cdot)$ is a solution to $\big( \mathcal{S}_{\left( u^{\prime
},b^{\prime }\right) }^{\prime }\big) $ while $\dot{\gamma}(t)> 0$ for a.e.\ $t\in \left[ 0,T\right] $, we get by using (\ref{ctrafo}) that 
\begin{align*}
-\dot{z}\left( t\right)=\dot{y}\left( \gamma (t)\right) \dot{\gamma}(t)&\in 
\dot{\gamma}(t)N_{C_{\left( u^{\prime },b^{\prime }\right) }(\gamma
(t))}\left( y(\gamma (t))\right)=N_{C_{\left( u^{\prime },b^{\prime
}\right) }(\gamma (t))}\left( z(t)\right)\\
&=N_{C_{\left( u,b\right)}(t)}\left( z(t)\right)\;\mathrm{ a.e.}\;t\in[0.T].
\end{align*}
It follows from \eqref{gamma} that $\gamma \in W^{1,1}\left( [0,T],\mathbb{R}\right) $, and so $z\in W^{1,1}([0,T],H)$ as well. Furthermore,  we have that $z\left( 0\right) =y\left( \gamma
(0)\right) =y(0)=x_{0}$ because $y(\cdot)$ is a solution of $\big( \mathcal{S}_{\left( u^{\prime },b^{\prime }\right) }^{\prime }\big)$. This allows us to conclude that $z(\cdot)$ is a solution of the original sweeping process $\big( \mathcal{S}_{\left( u,b\right)}\big) $ and---being absolutely continuous on $[0,T]$---it is unique by \cite[Theorem~3]{kunzemonteiro}. \qed
\end{proof}\vspace*{-0.05in}

Finally in this section, we present a consequence of Theorem~\ref{existsobolev} ensuring the result of this type for the $\dd-$ moving polyhedron \eqref{delmov}. This result is important to our applications to stability in the next section.\vspace*{-0.05in}
\begin{corollary}
\label{Cdelta} Let $\mathcal{H}$ be a separable Hilbert space, and let the uniform Slater condition \eqref{unifslater} be satisfied along a given control
$\left( \bar{u},\bar{b}\right) \in W^{1,1}([0,T],\mathcal{
H}^{m})\times W^{1,1}([0,T],\mathbb{R}^{m})$. Then there exists $\varepsilon >0$ such
that for all numbers $\delta \in \lbrack 0,\varepsilon )$ the perturbed sweeping process 
\begin{equation}
-\dot{x}\in N(C_{\left( \bar{u},\bar{b}\right) }^{\left( \delta \right)
}(t),x(t))\quad \mathrm{a.e.\;}t\in \left[ 0,T\right] ,\quad x(0)=\widehat{x}%
\left( 0\right) \in C_{\left( \bar{u},\bar{b}\right)}^{\left( \delta
\right) }(0)  \label{ydelta}
\end{equation}
admits a unique absolutely continuous solution. Here $C_{\left( \bar{u},\bar{b}
\right) }^{\left( \delta \right) }$ is defined in \eqref{delmov} and $
\widehat{x}(\cdot)$ is the continuous selection $\widehat{x}(t)\in C_{\left( \bar{u},\bar{b}\right) }^{\left( \delta \right) }(t)$ taken from \eqref{select1}.
\end{corollary}\vspace*{-0.15in}
\begin{proof}
As in the proof of Lemma~\ref{strongselection}, choose $\varepsilon >0$ from 
\eqref{slater2} and pick $\delta \in \lbrack 0,\varepsilon )$. Then $C_{\left( 
\bar{u},\widetilde{b}\right) }=C_{\left( \bar{u},\bar{b}\right) }^{\left(
\delta \right) }$, with $\widetilde{b}$ defined by $\widetilde{b}
_{i}:=b_{i}-\delta $ as $i=1,\ldots,m$, also satisfies the uniform Slater
condition. The result now follows from Theorem~\ref{existsobolev}. \qed
\end{proof}\vspace*{-0.25in}

\section{Quantitative stability of the perturbed sweeping dynamics}\label{quantstab}\vspace*{-0.05in}

\noindent In this section, we investigate the stability of solutions to
controlled polyhedral sweeping processes with respect to perturbations of controls and initial values of the sweeping dynamics. Theorem~\ref{existsobolev} allows
us to associate with each absolutely continuous control $\left( u,b\right)$
satisfying (\ref{unifslater}) and with the initial value $ x(0)=x_{0}\in C_{\left( u,b\right) }\left( 0\right)$ the unique absolutely
continuous solution $x_{(u,b)}$ of the sweeping process $\big( 
\mathcal{S}_{\left( u,b\right) }\big)$. In contrast with the previous
analysis, where the initial point $x_0$ was fixed, we now compare
solutions of $\big( \mathcal{S}_{\left( u,b\right) }\big)$ corresponding not
only to different controls but also to different initial points. To emphasize this dependence, let us write $\big(\mathcal{S}_{\left(
u,b,x_{0}\right) }\big)$ for the sweeping process $\big(\mathcal{S}
_{\left( u,b\right) }\big) $ corresponding to the initial condition $x(0)=x_{0}\in
C_{\left( u,b\right) }\left( 0\right) $ and denote its unique solution by $
x_{\left( u,b,x_{0}\right)}$. We begin with the following estimate,
which is based on Lemma~\ref{strongselection} and uses the arguments from the proof of Proposition~3 in \cite{HaddadJouraniThibault}.\vspace*{-0.05in} 
\begin{lemma}\label{LemmaEstim} Assume that $\mathcal{H}$ is a separable
Hilbert space, and that the uniform Slater condition \eqref{unifslater} holds for some given control $(\bar{u},\bar{b})\in W^{1,1}([0,T],\mathcal{H}
^{m})\times W^{1,1}([0,T],\mathbb{R}^{m})$. Then there exists $
\varepsilon>0$ such that for all $\delta \in \left( 0,\varepsilon \right) $,
for all controls $(u,b)\in \mathbb{B}_{1,1}\left( (\bar{u},\bar{b}),\frac{
\delta }{1+\Vert y_{\delta }\Vert_{\infty }}\right) $, and for all corresponding solutions $x(\cdot)$ to the sweeping processes $\big( \mathcal{S}_{\left(
u,b,x_{0}\right) }\big) $ we have the estimate
\begin{align}
&\Vert \dot{x}(t)\Vert \leq \frac{1}{\delta }\left( \Vert \widehat{x}\Vert
_{\infty }+\Vert y_{\delta }\Vert _{\infty }+\alpha _{\delta }\right) \left(
1+\Vert y_{\delta }\Vert _{\infty }+\alpha _{\delta }\right)
\sum_{i=1}^{m}\left( \Vert \dot{u}_{i}(t)\Vert +|\dot{b}_{i}(t)|\right)\notag\\
&\mathrm{a.e.}\,\,t\in \lbrack 0,T].  \label{EstTraj2}
\end{align}
Here $\widehat{x}(\cdot)$ stands for the continuous selection $\widehat{x}(t)\in C_{\left(\bar{u},\bar{b}\right) }^{\left( \delta \right) }(t)$ taken from \eqref{select1}, $y_{\delta}(\cdot)$ refers to the associate unique solution of the perturbed sweeping process \eqref{ydelta} guaranteed by Corollary~{\rm\ref{Cdelta}}, and the constant $\al_\dd$ is defined by
\begin{equation}
\alpha _{\delta }:=\int_{0}^{T}\Vert \dot{y}_{\delta }(t)\Vert dt+\sqrt{
\left( \int_{0}^{T}\Vert \dot{y}_{\delta }(t)\Vert dt\right) ^{2}+\Vert x(0)-
\widehat{x}\left( 0\right) \Vert ^{2}}.\label{alphadelta2}
\end{equation}
\end{lemma}\vspace*{-0.05in}
\begin{proof} As in previous proofs, we choose $\varepsilon >0$ from perturbed uniform Slater condition (\ref {slater2}) equivalent to the the assumed one
(\ref{unifslater}) by Proposition~\ref{slaterequiv}. Fix an arbitrary $\delta \in \left( 0,\varepsilon \right)$, then fix an arbitrary control pair
\begin{equation}
(u,b)\in \mathbb{B}_{1,1}\left( (\bar{u},\bar{b}),\frac{\delta }{1+\Vert
y_{\delta }\Vert _{\infty }}\right),  \label{newradius}
\end{equation}
and denote by $x(\cdot)$ the corresponding unique solution of the sweeping process $\big( 
\mathcal{S}_{\left( u,b,x_{0}\right) }\big)$ due to Theorem~\ref{existsobolev}. 
By the absolute continuity
of the triple $(u,b,x)$, the derivatives $\dot{x}(t)$, $\dot{u}_{i}(t)$ and $\dot{b}_{i}(t)$ exist for almost all $t\in \left[ 0,1\right] $. Fixing now any such time $t$ and then get 
\begin{align*}
x(t-s)&=x(t)-s(\dot{x}(t)+\alpha _{x}(s)),\quad u_{i}(t-s)=u_{i}(t)-s(\dot{u}%
_{i}(t)+\alpha _{u,i}(s))\\b_{i}(t-s)&=b_{i}(t)-s(\dot{b}_{i}(t)+\alpha
_{b,i}(s)),
\end{align*}
where $\lim_{s\rightarrow 0}\alpha _{x}(s)=0$, $\lim_{s\rightarrow 0}\alpha
_{u,i}(s)=0$ and $\lim_{s\rightarrow 0}\alpha _{b,i}(s)=0$. Since $x(t-s)\in
C_{\left( u,b\right) }\left( t-s\right) $ for all $s$, we deduce from (\ref
{select3}) that 
\begin{align*}
x(t-s)\in C_{(u,b)}(t)+\hskip 9cm&\\ \frac{1}{\delta }\Vert x(t-s)-\widehat{x}(t)\Vert
\sum_{i=1}^{m}\left( \Vert u_{i}(t-s)-u_{i}(t)\Vert \cdot \Vert x(t-s)\Vert
+|b_{i}(t-s)-b_{i}(t)\right) |\mathbb{B}&,
\end{align*}
where $\mathbb{B}$ refers as usual to the unit ball in $\mathcal{H}$. Using the
convexity of the $C_{(u,b)}(t)$ and passing to the limit $s\downarrow 0$, gives us the inclusion
\begin{equation*}
-\dot{x}(t)\in T(C_{(u,b)}(t),x(t))+\frac{1}{\delta }\Vert x(t)-\widehat{x}
(t)\Vert \sum_{i=1}^{m}\left( \Vert \dot{u}_{i}(t)\Vert \cdot \Vert
x(t)\Vert +|\dot{b}_{i}(t)|\right)\mathbb{B},
\end{equation*}
where $T(S,u)$ stands for the tangent cone to a convex set $S$ at $u$ in the sense of convex analysis. As $-\dot{x}(t)\in N(C_{(u,b)}(t),x(t))$, we arrive at
\begin{equation*}
\Vert \dot{x}(t)\Vert ^{2}\leq \Vert \dot{x}(t)\Vert \cdot \frac{1}{\delta }
\Vert x(t)-\widehat{x}(t)\Vert \sum_{i=1}^{m}\left( \Vert \dot{u}
_{i}(t)\Vert \cdot \Vert x(t)\Vert +|\dot{b}_{i}(t)|\right),
\end{equation*}
which in turn implies, since $t$ was arbitrarily chosen from a subset of
full measure on $[0,T]$, the derivative norm estimate
\begin{equation}
\Vert \dot{x}(t)\Vert \leq \\\frac{1}{\delta }\Vert x(t)-\widehat{x}\left(
t\right) \Vert \sum_{i=1}^{m}\left( \Vert \dot{u}_{i}(t)\Vert \cdot \Vert
x(t)\Vert +|\dot{b}_{i}(t)|\right) \quad \mathrm{a.e.\,}t\in \lbrack 0,T].
\label{EstTraj}
\end{equation}
To proceed further, let $y_{\delta}(\cdot)$ be the unique absolutely continuous solution to the sweeping process (\ref{ydelta}) according to Corollary~\ref{Cdelta}. Since $\left\langle \bar{u}_{i}\left( t\right) ,y_{\delta }(t)\right\rangle \le\bar{b}_{i}\left( t\right) -\delta $ for all $t\in \lbrack 0,T]$ and all $
i=1,\ldots ,m$, we deduce from (\ref{newradius}) that 
\begin{eqnarray*}
\langle u_{i}(t),y_{\delta }(t)\rangle -b_{i}(t) &\leq &\langle u_{i}(t)-%
\bar{u}_{i}(t),y_{\delta }(t)\rangle +\bar{b}_{i}(t)-b_{i}(t)-\delta \\
&\leq &\Vert u-\bar{u}\Vert _{\infty }\Vert y_{\delta }\Vert _{\infty
}+\Vert b-\bar{b}\Vert _{\infty }-\delta\\& \leq& \Vert u-\bar{u}\Vert
_{1,1}\Vert y_{\delta }\Vert _{\infty }+\Vert b-\bar{b}\Vert _{1,1}-\delta \\
&\leq &\Vert (u,b)-(\bar{u},\bar{b})\Vert _{1,1}(1+\Vert y_{\delta }\Vert
_{\infty })-\delta \leq 0\quad \forall t\in \lbrack 0,T].
\end{eqnarray*}
Therefore, $y_{\delta }(t)\in C_{(u,b)}(t)$ for all $t\in \lbrack 0,T]$. Remembering that $x(\cdot)$ solves the original sweeping process $\big( \mathcal{S}_{\left(
u,b,x_{0}\right) }\big) $, it follows that $-\dot{x}(t)\in N_{C_{\left( u,b\right) }(t)}\left(
x(t)\right) $ for a.e.\ $\mathrm{\,} t\in \lbrack 0,T]$, and hence we have
\begin{align*}
\frac{d}{dt}\frac{1}{2}\Vert x(t)-y_{\delta }(t)\Vert ^{2}& =
\langle \dot{x}
(t)-\dot{y}_{\delta }(t),x(t)-y_{\delta }(t)\rangle\\& =\langle \dot{x}
(t),x(t)-y_{\delta }(t)\rangle +\langle -\dot{y}_{\delta }(t),x(t)-y_{\delta
}(t)\rangle \\
& \leq \langle -\dot{y}_{\delta }(t),x(t)-y_{\delta }(t)\rangle \leq \Vert 
\dot{y}_{\delta }(t)\Vert \cdot \Vert x(t)-y_{\delta }(t)\Vert _{\infty }.
\end{align*}
This brings us to the estimate 
\begin{equation*}
\frac{\Vert x(t)-y_{\delta }(t)\Vert ^{2}}{2}-\frac{\Vert x(0)-\widehat{x}
\left( 0\right) \Vert ^{2}}{2}\leq \Vert x-y_{\delta }\Vert _{\infty }\cdot
\int_{0}^{T}\Vert \dot{y}_{\delta }(t)\Vert dt\quad \forall t\in \lbrack
0,T],
\end{equation*}
which implies on turn that
\begin{equation*}
\frac{\Vert x-y_{\delta }\Vert _{\infty }^{2}}{2}-\frac{\Vert x(0)-\widehat{x
}\left( 0\right) \Vert ^{2}}{2}\leq \Vert x-y_{\delta }\Vert _{\infty }\cdot
\int_{0}^{T}\Vert \dot{y}_{\delta }(t)\Vert dt.
\end{equation*}
Consequently, we arrive at the inequality
\begin{equation*}
\Vert x-y_{\delta }\Vert _{\infty }^{2}-2\left( \int_{0}^{T}\Vert \dot{y}
_{\delta }(t)\Vert dt\right) \Vert x-y_{\delta }\Vert _{\infty }-\Vert x(0)-
\widehat{x}\left( 0\right) \Vert ^{2}\leq 0.
\end{equation*}
Invoking the definition of $\al_\dd$ in \eqref{alphadelta2} gives us the estimate
\begin{equation}
\Vert x-y_{\delta }\Vert _{\infty }\leq \alpha _{\delta },
\label{alphadelta}
\end{equation}
which being combined with (\ref{EstTraj}) verifies the claimed
inequality (\ref{EstTraj2}) and thus completes the proof of the lemma. \qed
\end{proof}\vspace*{-0.05in}

Now we are ready to establish the main stability result.\vspace*{-0.07in}
\begin{theorem}
\label{controltostate} Let $\mathcal{H}$ be a separable Hilbert space, and let the uniform Slater condition \eqref{unifslater} hold for a given control pair
$(\bar{u},\bar{b})\in W^{1,1}([0,T],\mathcal{H}^{m})\times
W^{1,1}([0,T],\mathbb{R}^{m})$. Then there exist a number $\rho >0$ and a continuous function $K:\mathcal{H}\times\mathcal{H}\to\mathbb{R}_+$ such that for all control pairs
\begin{equation}
(u,b),(u^{\prime },b^{\prime })\in \left[ W^{1,1}([0,T],\mathcal{H}%
^{m})\times W^{1,1}([0,T],\mathbb{R}^{m})\right] \cap \mathbb{B}_{1,1}\left(
(\bar{u},\bar{b}),\rho \right) ,  \label{control ball}
\end{equation}
for all initial values $x_{0}\in C_{(u,b)}(0)$, $x_{0}^{\prime }\in
C_{(u^{\prime },b^{\prime })}(0)$, and the associated solutions $x,x^{\prime
}$ to the sweeping processes $\big( \mathcal{S}_{\left( u,b,x_{0}\right)
}\big)$ and $\big( \mathcal{S}_{\left( u^{\prime },b^{\prime },x_{0}^{\prime
}\right) }\big) $, respectively, we have\vspace*{-0.1in}
\begin{equation}
\left\Vert x(t)-x^{\prime }(t)\right\Vert ^{2}\leq \left\Vert
x_{0}-x_{0}^{\prime }\right\Vert ^{2}+K(x_0,x^\prime_0)\Vert (u-u^{\prime },b-b^{\prime
})\Vert _{\infty }\quad \forall t\in \lbrack 0,T].  \label{hoelder}
\end{equation}
\end{theorem}\vspace*{-0.15in}
\begin{proof} As above, we employ the equivalent description (\ref{slater2}) of the uniform Slater condition (\ref{unifslater}) and take $\ve>0$ from Proposition~\ref{slaterequiv}. Fixing an arbitrary number $\delta \in \left( 0,\varepsilon\right)$, define the quantity
\begin{equation}
\rho :=\min \left\{ \frac{\delta }{1+\Vert y_{\delta }\Vert _{\infty }},
\frac{\varepsilon -\delta }{3\left( 1+\left\Vert \widehat{x}\right\Vert
_{\infty }\right) }\right\},  \label{rhomindef}
\end{equation}
where $\widehat{x}(\cdot)$ is the continuous selection $\widehat{x}(t)\in C_{\left( \bar{u},\bar{b}\right) }^{\left( \delta \right) }(t)$ satisfying (\ref{select1}), and where $y_{\delta}(\cdot)$ is the unique absolutely continuous solution to the perturbed sweeping process (\ref{ydelta}) taken from
Corollary~\ref{Cdelta}. Select
arbitrary controls $(u,b),(u^{\prime },b^{\prime })$ from (\ref{control ball}), arbitrary initial values $x_{0}\in C_{(u,b)}(0)$, $x_{0}^{\prime }\in
C_{(u^{\prime },b^{\prime })}(0)$,  and the associated solutions $x,x^{\prime }
$ to the sweeping processes $\big(\mathcal{S}_{\left( u,b,x_{0}\right)
}\big)$ and $\big( \mathcal{S}_{\left( u^{\prime },b^{\prime },x_{0}^{\prime
}\right)}\big)$, respectively. Then it follows from (\ref{alphadelta}) that
\begin{equation}
\Vert x-y_{\delta }\Vert _{\infty }\leq \alpha _{\delta }\;\mbox{ and }\;\Vert
x^{\prime }-y_{\delta }\Vert _{\infty }\leq \alpha _{\delta }^{\prime }
\label{yaldel}
\end{equation}
for $\alpha _{\delta }$ defined in (\ref{alphadelta2}) and $\alpha
_{\delta }^{\prime }$ defined by the same formula with  the initial value $
x\left( 0\right) =x_{0}$ replaced by the initial value $x^{\prime }\left(
0\right) =x_{0}^{\prime }$.                                                                                                                                                                                                                                                                                                                 Lemma~\ref{LemmaEstim} gives us estimate (\ref{EstTraj2}) for the control $(u,b)$ as well as the corresponding estimate 
\begin{align}
&\Vert \dot{x}^{\prime }(t)\Vert \leq \delta ^{-1}\left( \Vert \widehat{x}%
\Vert _{\infty }+\Vert y_{\delta }\Vert _{\infty }+\alpha _{\delta }^{\prime
}\right) \left( 1+\Vert y_{\delta }\Vert _{\infty }+\alpha _{\delta
}^{\prime }\right) \sum_{i=1}^{m}\left( \Vert \dot{u}_{i}^{\prime }(t)\Vert
+|\dot{b}_{i}^{\prime }(t)|\right)\notag\\
&\mathrm{a.e.}\,\,t\in \lbrack 0,T]\label{EstTraj6}
\end{align}
for the control $(u^{\prime },b^{\prime })$. Denoting now
\begin{align}
C:=&\left( \alpha _{\delta }+\Vert y_{\delta }\Vert _{\infty }+\Vert {
\widehat{x}}\Vert _{\infty }\right) \left( 1+\alpha _{\delta }+\Vert
y_{\delta }\Vert _{\infty }\right), \notag\\C^{\prime }:=&\left( \alpha
_{\delta }^{\prime }+\Vert y_{\delta }\Vert _{\infty }+\Vert {\widehat{x}}
\Vert _{\infty }\right) \left( 1+\alpha _{\delta }^{\prime }+\Vert y_{\delta
}\Vert _{\infty }\right)  \label{ccbar}
\end{align}
and integrating  (\ref{EstTraj2}) ensure that 
\begin{equation*}
\int_{0}^{t}\Vert \dot{x}(s)\Vert ds\leq \delta
^{-1}C\sum_{i=1}^{m}\int_{0}^{t}\left( \Vert \dot{u}_{i}(s)\Vert +|\dot{b}
_{i}(s)|\right) ds\quad \forall \,t\in \lbrack 0,T].
\end{equation*}
Therefore, recalling that $(u,b)\in \mathbb{B}_{1,1}\left( (\bar{u},\bar{b}
),\rho \right)$ yields
\begin{equation}
\int_{0}^{t}\Vert \dot{x}(s)\Vert ds\leq \delta ^{-1}C\left\Vert \left(
u,b\right) \right\Vert _{1,1}\leq \delta ^{-1}C\left( \rho +\left\Vert
\left( \bar{u},\bar{b}\right) \right\Vert _{1,1}\right) .  \label{EstTraj3}
\end{equation}
Similarly, the integration of (\ref{EstTraj6}) gives us
\begin{equation}
\int_{0}^{t}\Vert \dot{x}^{\prime }(s)\Vert ds\leq \delta ^{-1}C^{\prime
}\left( \rho +\left\Vert \left( \bar{u},\bar{b}\right) \right\Vert
_{1,1}\right). \label{EstTraj4}
\end{equation}
Let now $t\in \lbrack 0,T]$ be from a subset of full measure such that $\dot{
x}(t)$ and $\dot{x}^{\prime }(t)$ exist. We clearly have $x(t)\in C_{(u,b)}(t)$ and 
$x^{\prime }(t)\in C_{(u^{\prime },b^{\prime })}(t)$. Since a ball in the $
\Vert \cdot \Vert _{1,1}$-norm is contained in a ball of the same radius in
the $\Vert \cdot \Vert _{\infty }$-norm, the construction of $\rho $ in (\ref
{rhomindef}) allows us to employ the error bound (\ref{select3}) from Lemma \ref
{strongselection}. This ensures the existence of $x_{1}\in C_{(u^{\prime
},b^{\prime })}(t)$ and $x_{1}^{\prime }\in C_{(u,b)}(t)$ with 
\begin{equation*}
\begin{array}{ll}
\Vert x(t)-x_{1}\Vert\leq\delta ^{-1}{\Vert x(t)-\widehat{x}(t)\Vert }
\disp\max_{i=1,\ldots ,m}\left[ \langle u_{i}^{\prime }(t)-u_{i}(t),x\rangle
+b_{i}(t)-b_{i}^{\prime }(t)\right] _{+}\\
\leq\delta ^{-1}{\Vert x(t)-\widehat{x}(t)\Vert }\left( \Vert
u(t)-u^{\prime }(t)\Vert \left\Vert x{(t)}\right\Vert +\left\Vert
b(t)-b^{\prime }(t)\right\Vert \right)\\
\leq\delta ^{-1}{\Vert x(t)-\widehat{x}(t)\Vert }(1+\Vert x(t)\Vert
)(\Vert (u(t)-u^{\prime }(t),b(t)-b^{\prime }(t))\Vert.
\end{array}
\end{equation*}
Similar considerations bring us to the estimate
\begin{equation*}
\Vert x^{\prime }(t)-x_{1}^{\prime }\Vert \leq \delta ^{-1}{\Vert x}^{\prime
}{(t)-\widehat{x}(t)\Vert }(1+\Vert x^{\prime }(t)\Vert )(\Vert
(u(t)-u^{\prime }(t),b(t)-b^{\prime }(t))\Vert.
\end{equation*}
Since $x(\cdot)$ and $x^{\prime }(\cdot)$ are absolutely continuous solutions to $\big( \mathcal{S}_{\left(u,b,x_{0}\right)}\big)$ and
$\big(\mathcal{S}_{\left( u^{\prime },b^{\prime
},x_{0}^{\prime }\right) }\big)$, respectively, we deduce from $-\dot{x}
(t)\in N_{C_{\left( u,b\right) }(t)}\left( x(t)\right) $, $-\dot{x}
^{\prime }(t)\in N_{C_{\left( u^{\prime },b^{\prime }\right) }(t)}\left(
x^{\prime }(t)\right) $, and the obtained estimates of $\Vert x(t)-x_{1}\Vert$ and $\Vert x^{\prime }(t)-x_{1}^{\prime }\Vert$ that
\begin{eqnarray*}
&&{\frac{d}{{dt}}}\frac{1}{2}\left\Vert x(t)-x^{\prime }(t)\right\Vert ^{2}
\\
&=&\langle \dot{x}(t)-\dot{x}^{\prime }(t),x(t)-x^{\prime }(t)\rangle
=\langle \dot{x}(t),x(t)-x^{\prime }(t)\rangle -\langle \dot{x}^{\prime
}(t),x(t)-x^{\prime }(t)\rangle  \\
&=&\langle \dot{x}(t),x(t)-x_{1}^{\prime }\rangle +\langle \dot{x}
(t),x_{1}^{\prime }-x^{\prime }(t)\rangle +\langle \dot{x}^{\prime
}(t),x^{\prime }(t)-x_{1}\rangle +\langle \dot{x}(t),x_{1}-x(t)\rangle  \\
&\leq &\langle \dot{x}(t),x_{1}^{\prime }-x^{\prime }(t)\rangle +\langle 
\dot{x}(t),x_{1}-x(t)\rangle\\ &\leq& \Vert \dot{x}(t)\Vert \Vert x_{1}^{\prime
}-x^{\prime }(t)\Vert +\Vert \dot{x}^{\prime }(t)\Vert \Vert x_{1}-x(t)\Vert 
\\
&\leq &\delta ^{-1}\left( \Vert \dot{x}(t)\Vert {\Vert x^{\prime }(t)-
\widehat{x}(t)\Vert }(1+\Vert x^{\prime }(t)\Vert )+\Vert \dot{x}^{\prime
}(t)\Vert {\Vert x(t)-\widehat{x}(t)\Vert }(1+\Vert x(t)\Vert )\right)  \\
&&\cdot \Vert (u(t)-u^{\prime }(t),b(t)-b^{\prime }(t))\Vert.
\end{eqnarray*}
For all $t\in[0,T]$ define the function
\begin{equation*}
\chi \left( t\right) :=\delta ^{-1}\left( {\Vert x^{\prime }(t)-\widehat{x}
(t)\Vert }(1+\Vert x^{\prime }(t)\Vert )+{\Vert x(t)-\widehat{x}(t)\Vert }
(1+\Vert x(t)\Vert )\right).
\end{equation*}
Then the latter estimate can be rewritten as 
\begin{equation}
{\frac{d}{{dt}}}\frac{1}{2}\left\Vert x(t)-x^{\prime }(t)\right\Vert
^{2}\leq \chi \left( t\right) \left( \Vert \dot{x}(t)\Vert +\Vert \dot{x}
^{\prime }(t)\Vert \right) \Vert (u(t)-u^{\prime }(t),b(t)-b^{\prime
}(t))\Vert.  \label{chiineq}
\end{equation}
It follows from (\ref{yaldel}) and (\ref{ccbar}) that $\chi \left( t\right)
\leq \delta ^{-1}\left( C+C^{\prime }\right) $.
As $t$ was arbitrarily chosen from a subset of full measure of $\left[ 0,T
\right]$, we integrate (\ref{chiineq}) and then employ (\ref{EstTraj3})
and (\ref{EstTraj4}) to get
\begin{eqnarray*}
&&\left\Vert x(t)-x^{\prime }(t)\right\Vert ^{2}-\left\Vert x(0)-x^{\prime
}(0)\right\Vert ^{2}\\ &\leq &\delta ^{-1}\left( C+C^{\prime }\right)
\int_{0}^{t}\left( \Vert \dot{x}(s)\Vert +\Vert \dot{x}^{\prime }(s)\Vert
\right) \Vert (u(s)-u^{\prime }(s),b(s)-b^{\prime }(s))\Vert ds \\
&\leq &\delta ^{-1}\left( C+C^{\prime }\right) \Vert (u-u^{\prime
},b-b^{\prime })\Vert _{\infty }\int_{0}^{t}\left( \Vert \dot{x}(s)\Vert
+\Vert \dot{x}^{\prime }(s)\Vert \right) ds \\
&\leq &\underbrace{\delta ^{-2}\left( C+C^{\prime }\right) ^{2}\left( \rho +\left\Vert
\left( \bar{u},\bar{b}\right) \right\Vert _{1,1}\right)}_{K(x_0,x^\prime_0)} \Vert (u-u^{\prime},b-b^{\prime })\Vert _{\infty }
\end{eqnarray*}
for all $t\in \lbrack 0,T]$. By (\ref{ccbar}), $C$ and $C^\prime$ depend continuously on $\alpha_\delta$ and $\alpha^\prime_\delta$, respectively, which in turn depend continuously on $x_0$ and $x^\prime_0$ by (\ref{alphadelta2})).
Thus we verify that the obtained continuous function $K(x_0,x\prime_0)$
ensures the claimed estimate \eqref{hoelder}, and we are done with the proof of the theorem. \qed
\end{proof}\vspace*{-0.05in}

To conclude this section, we present a direct consequence of 
Theorem~\ref{controltostate} for the case where the initial value $x_0$ in \eqref{sweeping} is fixed. In this case the function $K(\cdot)$ in the estimate
\eqref{hoelder} is constant.\vspace*{-0.05in}  
\begin{corollary}
Let $\mathcal{H}$ be a separable Hilbert space, let the uniform Slater condition \eqref{unifslater} hold  for a given control
$(\bar{u},\bar{b})\in W^{1,1}([0,T],\mathcal{H}^{m})\times
W^{1,1}([0,T],\mathbb{R}^{m})$, and let $x_0\in C_{(\bar{u},\bar{b})}$ be an arbitrarily given initial value in \eqref{sweeping}. Then there exist positive numbers $\rho$ and $K$ such that for all controls $(u,b),(u^{\prime },b^{\prime })$ satisfying \eqref{control ball} and the corresponding solutions $x(\cdot)$ and $x^{\prime
}(\cdot)$ of the sweeping processes $\big(\mathcal{S}_{\left( u,b,x_{0}\right)
}\big)$ and $\big( \mathcal{S}_{\left( u^{\prime },b^{\prime },x_{0}\right) }\big) $ with $x_{0}\in C_{(u,b)}(0)\cap C_{(u^{\prime },b^{\prime })}(0)$, respectively, we have
\begin{equation*}
\left\Vert x(t)-x^{\prime }(t)\right\Vert ^{2}\leq K\Vert (u-u^{\prime },b-b^{\prime
})\Vert _{\infty }\quad \forall t\in \lbrack 0,T].  
\end{equation*}
\end{corollary}\vspace*{-0.25in}

\section{Discrete approximations of controlled sweeping
processes}\label{discapp} \setcounter{equation}{0}\vspace*{-0.05in}

\noindent The last two sections of the paper are devoted to the study of the
following \emph{optimal control problem} for the sweeping process \eqref{sweeping} with controls in polyhedral moving sets \eqref{movpoly} and additional \emph{endpoint constraints} as well as the {\em pointwise equality constraints} on the \emph{$u$-control} functions:
\begin{align}
\min \left\{ f(u,b)|(u,b)\in W^{1,2}([0,T],\mathbb{R}^{nm}\times \mathbb{R}
^{m}),\,\,\Vert u_{i}(t)\Vert =1\,\,(i=1,\ldots ,m)\right.&\notag\\ \left. x_{\left(
u,b\right) }(T)\in \Omega \right\},&  \tag{P}
\end{align}
where $\Omega \subseteq\mathbb{R}^{n}$ is a closed subset, $f$ is a
cost function (specified later on), and $x_{\left(u,b\right) }$ is
the unique trajectory of the polyhedral sweeping process
$\big(\mathcal{S}_{(u,b)}\big)$ from \eqref{sweeping} generated by a control
pair $(u,b)=(u(\cdot),b(\cdot))$ of the above class on $[0,T]$. Such
a control pair $(u,b)$ is called a \emph{feasible solution} to ($P$) if $\|u(t)\|$=1 for all $t\in[0,T]$ and $x_{(u,b)}(T)\in\Omega$ for the corresponding trajectory of \eqref{sweeping}. Note that our focus in what follows is on \emph{Lipschitzian} controls in $(P)$, which uniquely generate by Theorem~\ref{existlip} Lipschitzian sweeping trajectory under the imposed \emph{uniform Slater condition} \eqref{unifslater}. At the current stage of our developments for $(P)$, we have to restrict ourselves to the case of finite-dimensional state spaces.

Our main goal here is to develop the {\em method of discrete
approximations} to investigate the sweeping control problem $(P)$ and its
discrete counterparts from both viewpoints of {\em stability} and
deriving {\em necessary suboptimality} and {\em optimality
conditions}. Stability issues address the construction of
finite-difference approximations of sweeping differential inclusions
such that their feasible solutions strongly approximate a broad
class of {\em canonical controls} in the original sweeping process;
this notion is introduced in the paper for the {\em first time}.
Furthermore, we construct a sequence of discrete-time optimal
control problems $(P_k)$ always admitting optimal solutions, which
$W^{1,2}$-{\em strongly converge} to a prescribed {\em local
minimizer} of the {\em intermediate} class (between weak and strong,
including the latter) of the original sweeping control problem
$(P)$. This opens the door to derive {\em necessary optimality
conditions} for such minimizers of $(P)$ by using advanced tools of
variational analysis and (first-order and second-order) generalized
differentiation. Furnishing this approach, we concentrate here on
deriving necessary optimality conditions for problems $(P_k)$ with
the approximation number $k\in\N$ being sufficiently large. The
obtained necessary optimality conditions for $(P_k)$ serve as
constructive {\em suboptimality} conditions for intermediate local
minimizers of $(P)$ that are convenient for numerical
implementations. This is a \emph{clear advantage} of the method of discrete
approximations in comparison with other methods of deriving
necessary optimality conditions for continuous-time variational and
control problems. In our separate publication, we are going to
realize the involved limiting procedure of passing to the limit from
the obtained necessary optimality conditions for $(P_k)$ (i.e.,
suboptimality conditions for $(P)$) to derive {\em exact} necessary
optimality conditions for intermediate local minimizers of
continuous-time sweeping control problems of type $(P)$.

The method of discrete approximations was developed in
\cite{m95,m06} to establish necessary suboptimality and optimality
conditions for {\em Lipschitzian} differential inclusions. Sweeping
differential inclusions are highly {\em discontinuous}, and the
machinery of Lipschitzian variational analysis is not applicable in
the sweeping framework. Further developments of this method in various
sweeping control settings can be found in
\cite{ao,mor1,cm,cg,mor2,cmn} and the references therein. However,
neither these publications, nor those of \cite{bk,pfs,zeidan} exploring other approaches to deriving optimality conditions in different models of sweeping optimal control address additional endpoint constraints $x(T)\in\O$ on sweeping trajectories.

In this section we focus on the construction of discrete
approximations for the constrained sweeping dynamics and local minimizers of $(P)$ with obtaining stability/convergence results, while the
next section is devoted to reviewing the required tools of generalized differentiation and their applications to necessary optimality conditions for
discrete approximation problems $(P_k)$ giving us suboptimality
conditions for intermediate local minimizers of
$(P)$.\vspace*{0.03in}

Let us start with introducing a new notion of {\em canonical
controls} for problem $(P)$ that plays a crucial role in our developments.\vspace*{-0.07in}

\begin{definition}\label{canon} We say that a control pair $(u,b)\in W^{1,2}([0,T],\mathbb{R}^{nm}\times\mathbb{R})$ is {\sc canonical} for problem $(P)$ if the
following conditions hold:\\
$\bullet$ The functions $u(\cdot)$ and $b(\cdot))$ are Lipschitz
continuous on $[0,T]$.\\
$\bullet$ The uniform Slater condition \eqref{unifslater} is satisfied along $(u,b)$.\\
$\bullet$ We have the constraints
\begin{equation*}
\|u_i(t)\|=1\;\mbox{ for all }\;t\in[0,T]\;\mbox{ and
}\;i=1,\ldots,m.
\end{equation*}
$\bullet$ The derivatives $\dot u(\cdot)$ and $\dot b(\cdot)$ are of
bounded variation $(BV)$ on $[0,T]$ together with the derivative of
the unique Lipschitz continuous trajectory $x(\cdot)$ of
\eqref{sweeping} generated by the control pair $(u,b)$.
\end{definition}\vspace*{-0.1in}

Observe that the corresponding trajectory to \eqref{sweeping}
generated by a canonical control pair {\em may not} satisfy the
endpoint constraint $x_{\left(u,b\right)}(T)\in\Omega$, i.e., not
any canonical pair is feasible for $(P)$.

To proceed with our approach, we construct a sequence of discrete
approximations of the sweeping process $(\mathcal{S}_{(u,b)})$ from
\eqref{sweeping} over the controlled polyhedron \eqref{movpoly} {\em
without any appeal to optimization} as in $(P)$. For each
$k\in\mathbb{N}$ define the discrete mesh on $[0,T]$ by
\begin{equation}  \label{e:DP}
\Delta_k:=\left\{0=t^k_0<t^k_1<\ldots<t^k_{\nu(k)-1}<t^k_{\nu(k)}=T\right\}
\end{equation}
with $h^k_j:=t^k_{j+1}-t^k_j\downarrow 0$, $j=0,\ldots,\nu(k)-1$, as 
$k\to\infty$. Denote
\begin{equation} \label{F}
F(u,b,x):=N_{C(u,b)},\;C(u,b):=\big\{x\in\mathbb{R}^n\;\big|
\;\left\langle u_i,x\right\rangle\le b_i\;\;(i=1,\ldots,m)\big\}.
\end{equation}

The following theorem tells us that \emph{any canonical} control
pair $(u,b)$ and the corresponding sweeping trajectory
$x(\cdot)$ can be \emph{$W^{1,2}$-strongly} approximated by feasible
solutions to discrete sweeping processes defined on the
partition $\Delta_k$ from \eqref{e:DP} and appropriately extended to
the continuous-time interval $[0,T]$.\vspace*{-0.05in}
\begin{theorem}
\label{da-feas} Let $\left(\bar{u}(\cdot),\bar b(\cdot)\right)$ be a
canonical control pair for $(P)$, and let $\ox(\cdot)$ be the
corresponding unique solution of the Cauchy problem in
\eqref{sweeping}. Then there exist a mesh $\Delta_k$ in
\eqref{e:DP}, a sequence of piecewise linear functions $(\widetilde
u^k(\cdot),\widetilde b^k(\cdot),\widetilde x^k(\cdot))$ on $[0,T]$,
and a sequence of positive numbers $\dd_k\downarrow 0$ as
$k\to\infty$ such that $(\widetilde u^k(0),\widetilde
b^k(0),\widetilde x^k(0))=(\bar{u}(0),\bar b(0),x_0)$,
\begin{equation}\label{e:a-dc}
1-\dd_k\le\left\|\widetilde u^k_i(t^k_j)\right\|\le 1+\dd_k\;%
\mbox{
for all }\;t^k_j\in\Delta_k,\quad i=1,\ldots,m,
\end{equation}
\begin{equation*}
\widetilde x^k(t)=\widetilde x^k(t^k_j)+(t-t^k_j)\widetilde
v^k_j,\;\;t^k_j\le t\le t^k_{j+1}\;\;\mbox{with}\;\;-\widetilde
v^k_j\in F\big(\widetilde u^k(t^k_j),\widetilde
b^k(t^k_j),\widetilde x^k(t^k_j)\big)
\end{equation*}
for $j=0,\ldots,\nu(k)-1$, $k\in\mathbb{N}$, and the sequence
$\{(\widetilde u^k(\cdot),\widetilde b^k(\cdot),\widetilde
x^k(\cdot))\}$ converges to $(\bar{u}(\cdot),\bar
b(\cdot),\bar{x}(\cdot))$ as $k\to\infty$ in the $W^{1,2}$-norm
topology on $[0,T]$.
\end{theorem}\vspace*{-0.2in}
\begin{proof}
As mentioned, the existence of the unique Lipschitz continuous
trajectory $\bar{x}(\cdot)$ of the Cauchy problem for the polyhedral
sweeping process in \eqref{sweeping} generated by the given
canonical control pair $(\bar{u}(\cdot),\bar b(\cdot))$ follows from
Theorem~\ref{existlip}. Now we are in a position to deduce the
claimed assertions from \cite[Theorem~4.1]{mor1} under the BV
assumption on $\dot{\ou}(\cdot)$, $\dot{\ob}(\cdot)$, and
$\dot{\ox}(\cdot)$. Indeed, the qualification condition (H4) from
\cite[Theorem~4.1]{mor1} is equivalent to the uniform Slater
condition \eqref{unifslater} by our new result obtained in
Proposition~\ref{slaterequiv}. Thus the application of
\cite[Theorem~4.1]{mor1} gives us all the assertions claimed in this
theorem. \qed
\end{proof}\vspace*{-0.1in}

From now on, we consider for simplicity problem $(P)$, where the cost
function is defined in the \emph{Mayer form} via a given terminal state
function $\varphi\colon\mathbb{R}^n\to\mathbb{R}$ by
\begin{equation*}
f(u,b):=\varphi\big(x_{u,b}(T)\big).
\end{equation*}
If the function $\varphi$ is lower semicontinuous, then problem
$(P)$
admits a (global) \emph{optimal solution} in $W^{1,2}([0,T],\mathbb{R}^{nm}\times%
\mathbb{R}^m)$ provided that there is a bounded minimizing sequence
of feasible solutions; see \cite[Theorem~3.1]{mor1} and its proof.
Since our main attention is paid to deriving necessary
(sub)optimality conditions in $(P)$, it is natural to define an
appropriate notion of \emph{local minimizers}.\vspace*{0.03in}

The notion of local minimizers of our study in this paper occupies an \emph{%
intermediate} position between the classical notions of \emph{weak} and
\emph{strong} minimizers in variational and control problems, while
encompassing the latter. Following \cite{m95}, where this notion was
initiated for Lipschitzian differential inclusions (see also \cite{m06} for
more details), we keep the name ``intermediate" for the version of this
notion in the setting of our sweeping control problem $(P)$.\vspace*{-0.1in}

\begin{definition}
\label{ilm} We say that a feasible control pair $(\bar{u},\bar b)$ for $(P)$
is an {\sc intermediate local minimizer} in this problem if there exists $%
\varepsilon>0$ such that
\begin{equation*}
\varphi\big(x_{\bar{u},\bar b}(T)\big)\le\varphi\big(x_{u,b}(T)\big)
\end{equation*}
for any feasible solution to $(P)$ satisfying the condition
\begin{equation} \label{loc}
\|(u,b)-(\bar{u},\bar b)\|_{W^{1,2}}+\|x_{u,b}-x_{\bar{u},\bar
b}\|_{W^{1,2}}\le\varepsilon.
\end{equation}
\end{definition}

The notion of \emph{strong local minimizer} for $(P)$ is a particular case
of Definition~\ref{ilm}, where the norm $\|x_{u,b}-x_{\bar{u},\bar
b}\|_{W^{1,2}}$ in \eqref{loc} is replaced with the norm $\|x_{u,b}-x_{\bar{u%
},\bar b}\|_{\mathcal{C}}$ in the space of continuous functions $\mathcal{C}%
([0,T],\mathbb{R}^n)$. It is not hard to construct examples showing
that there exist intermediate local minimizers to $(P)$ that fail to
be strong ones; see \cite{m95,m06,vinter} even for simpler
problems.\vspace*{0.05in}

Having $F(u,b,x)$ from \eqref{F}, fix a Lipschitz continuous
intermediate local minimizer $(\bar{u},\bar b)$ for $(P)$ with the
corresponding sweeping trajectory $\bar{x}(\cdot):=x_{\bar{u},\bar
b}$ and assume that the uniform Slater condition \eqref{unifslater}
holds along $(\bar{u},\bar b)$. Take the mesh $\Delta_k$ from
\eqref{e:DP} and identify the points $t^k_j$ with the index $j$ if
no confusion arises. Consider now discrete triples $(u^k,b^k,x^k)$
with the components
\begin{equation*}
(u^k,b^k,x^k):=(u^k_0,u^k_1,\ldots,u^k_{\nu(k)},b^k_0,b^k_1,\ldots,b^k_{%
\nu(k)},x^k_0,x^k_1,\ldots,x^k_{\nu(k)})
\end{equation*}
and form the sequence of \emph{discrete approximation} problems $(P_k)$ by: 
\begin{align} 
&\mbox{minimize}\quad\varphi\big(x^k_{\nu(k)})+\label{disc-cost}\\ &\frac{1}{2}\sum_{j=0}^{\nu(k)-1}
\int_{t^k_j}^{t^k_{j+1}}\bigg\|\bigg(\dfrac{u^k_{j+1}-u^k_j}{h^k_j},\dfrac{
b^k_{j+1}-b^k_j}{h^k_j},\dfrac{x^k_{j+1} -x^k_j}{h^k_j}\bigg)-\big(\dot{\bar{
u}}(t),\dot{\bar b}(t),\dot{\ox}(t)\big)\bigg\|^2dt\notag
\end{align}
over the triples $(u^k,b^k,x^k$) subject to the following constraints:
\begin{equation} \label{disc-sw}
x^k_{j+1}\in x^k_j-h^k_jF(u^k_j,b^k_j,x^k_j),\;j=0,\ldots,\nu(k)-1,
\end{equation}
\begin{equation}\label{ini}
x^k_0=x_0\in C_{\ou,\ob}(0),\;(u^k_0,b^k_0)=\big(\bar{u}(0),\bar b(0)\big)%
,\;x^k_{\nu(k)}\in\Omega+\xi_kI\!\!B,
\end{equation}
\begin{equation}\label{u-const}
1-\dd_k\le\left\|u^k_i(t^k_j)\right\|\le 1+\dd_k\;\mbox{ for all }
\;t^k_j\in\Delta_k,\quad i=1,\ldots,m,
\end{equation}
\begin{equation}  \label{ic1}
\sum_{j=0}^{\nu(k)-1}\int_{t^k_j}^{t^k_{j+1}}\Big\|\big(u^k_j,b^k_j,x^k_j
\big)-\big(\bar{u}(t),\bar
b(t),\bar{x}(t)\big)\Big\|^2dt\le\dfrac{\varepsilon}{2},
\end{equation}
\begin{equation}\label{ic2}
\sum_{j=0}^{\nu(k)-1}\int_{t^k_j}^{t^k_{j+1}}\bigg\|\bigg(\dfrac{%
u^k_{j+1}-u^k_j}{h^k_j},\dfrac{b^k_{j+1}-b^k_j}{h^k_j},\dfrac{x^k_{j+1}-x^k_j
}{h^k_j}\bigg)- \big(\dot{\bar{x}}(t),\dot{\bar a}(t),\dot{\bar b}(t)\big)
\bigg\|^2dt\le\dfrac{\varepsilon}{2},
\end{equation}
where $\{\dd_k\}$ in \eqref{u-const} is taken from
Theorem~\ref{da-feas} applied to $(\bar{u},\bar b)$ and can be chosen
such that both inequalities in \eqref{u-const} are strict, where
$\varepsilon>0$ in \eqref{ic1} and \eqref{ic2} is taken from
Definition~\ref{ilm} of the intermediate local minimizer
$(\bar{u},\bar b)$ for $(P)$, and where the sequence $\{\xi_k\}$ of
the endpoint perturbations in \eqref{ini} is defined by
\begin{equation} \label{end-pert}
\xi_k:=\|\widetilde x^k(T)-\bar{x}(T)\|\to 0\;\mbox{ as }\;k\in\mathbb{N}
\end{equation}
via the sequence $\{\widetilde x^k(\cdot)\}$ approximating $\bar{x}(\cdot)$
in Theorem~\ref{da-feas}.\vspace*{0.03in}

The next theorem establishes the existence of optimal solutions to
problems $(P_k)$ for all $k\in\mathbb{N}$ and then shows that any
sequence of optimal controls $\{(\bar{u}^k,\bar b^k)\}$ to $(P_k)$
constructed for the given {\em canonical intermediate local
minimizer} $(\ou,\ob)$ of $(P)$, together with the corresponding
sequence of discrete trajectories $\{\bar{x}^k\}$ piecewise linearly
extended to the whole interval $[0,T]$, \emph{strongly
$W^{1,2}$-converge} as $k\to\infty$ to the prescribed local optimal
triple $(\bar{u},\bar b,\bar{x})$ for $(P)$.\vspace*{-0.07in}

\begin{theorem}\label{ilm-conver} Let $(\bar{u},\bar b)$ be a canonical intermediate local minimizer
for $(P)$ with the corresponding sweeping trajectory
$\bar{x}(\cdot)$. The following assertions hold:\\[0.5ex]
{\bf(i)} If the cost function $\varphi$ is lower semicontinuous around $%
\bar{x}(T)$, then each problem $(P_k)$ admits an optimal solution whenever $%
k\in\mathbb{N}$ is sufficiently large.\newline {\bf(ii)} If in
addition $\varphi$ is continuous around $\bar{x}(T)$, then every
sequence of optimal solutions $\{(\bar{u}^k,\bar b^k)\}$ to $(P_k)$
and the corresponding sequence of discrete trajectories
$\{\bar{x}^k\}$,
being piecewise linearly extended to $[0,T]$, converge to $(\bar{u},\bar b,%
\bar{x})$ as $k\to\infty$ in the norm topology of $W^{1,2}([0,T],\mathbb{R}%
^{mn}\times\mathbb{R}^m\times\mathbb{R}^n)$.
\end{theorem}\vspace*{-0.2in}
\begin{proof}
To verify (i), observe first that the set of feasible solutions to
problem $ (P_k)$ is {\em nonempty} for all large $k\in\mathbb{N}$. Namely, we show that the approximating sequence
$\{(\widetilde u^k,\widetilde b^k,\widetilde x^k)\}$ from
Theorem~\ref{da-feas}, being applied to the given {\em canonical}
intermediate local minimizer $(\bar{u},\bar b)$ of the original
problem $(P)$, consists of feasible solutions to $(P_k)$ when $k$ is
sufficiently large. Indeed, the
discrete sweeping inclusions \eqref{disc-sw} with the initial data in %
\eqref{ini} are clearly satisfied for $\{(\widetilde u^k,\widetilde
b^k,\widetilde x^k)\}$ together with the control constraints
\eqref{u-const}, the conditions in \eqref{ic1} and \eqref{ic2} also
hold for large $k$ by the $W^{1,2}$-convergence of the extended
triples $\{(\widetilde u^k(t),\widetilde b^k(t),\widetilde
x^k(t))\}$ to $(\bar{u}(t),\bar b(t), \bar{x}(t))$ on $[0,T]$ as
$k\to\infty$, and the fulfillment of the endpoint constraint in
\eqref{ini} for the approximating trajectories $\widetilde
x^k(\cdot)$ follows from $\bar{x}(T)\in\O$ and the definition of
$\xi_k$ in \eqref{end-pert} by Theorem~\ref{da-feas} applied to the canonical intermediate local minimizer $(\ou,\ob)$. It follows from the construction of $(P_k)$ and the structure of $F$ in \eqref{F} that the set of
feasible solutions to $(P_k)$ is closed. Furthermore, the
constraints in \eqref{u-const}--\eqref{ic2} ensure the boundedness
of the latter set. Since $\varphi$ is assumed to be lower
semicontinuous around $\bar{x}(T)$, the existence of optimal
solutions to such $(P_k)$ follows from the classical Weierstrass
existence theorem in finite dimensions.

Now we verify assertion (ii) of the theorem. Consider an arbitrary sequence $\{(\bar{
u}^k(\cdot),\bar b^k(\cdot),\bar{x}^k(\cdot))\}$ of optimal controls to $
(P_k)$ and the associated trajectories of \eqref{disc-sw} that are piecewise
linearly extended to $[0,T]$. We aim at proving
\begin{equation} \label{lim-con}
\lim_{k\to\infty}\int^T_0\big\|(\dot{\ou}^k(t),\dot{\ob}^k(t),\dot{\ox}
^k(t)\big)-(\dot{\ou}(t),\dot{\ob}(t),\dot{\ox}(t)\big)\big\|^2dt=0,
\end{equation}
which readily yields the claimed convergence in (ii). Supposing on
the contrary that \eqref{lim-con} fails gives us a subsequence of $k\to\infty$ (no relabeling) along which the limit in 
\eqref{lim-con} equals to some $\sigma>0$. Due to \eqref{ic2}, the sequence $
\{(\dot{\ou}^k(t),\dot{\ob}^k(t),\dot{\ox}^k(t))\}$ is weakly
compact in $L^2([0,T],\mathbb{R}^{mn}\times\mathbb{R}^m\times\mathbb{R}^n)$,
and hence it contains a subsequence that converges to some triple $
(\vartheta^u(\cdot),\vartheta^b(\cdot),\vartheta^x(\cdot))\in L^2([0,T],
\mathbb{R}^{mn}\times\mathbb{R}^m\times\mathbb{R}^n)$ weakly in this
space. Employing Mazur's weak closure theorem tells us that there is a sequence of convex combinations of
$(\dot{\ou}^k(\cdot),\dot{\ob}^k(\cdot),\dot{\ox}^k(\cdot))$, which
converges to
$(\vartheta^u(\cdot),\vartheta^b(\cdot),\vartheta^x(\cdot))$ strongly in $L^2([0,T],\mathbb{R}^{mn}\times\mathbb{R}
^m\times\mathbb{R}^n)$, and hence almost everywhere on $[0,T]$ along
a subsequence. Define
\begin{equation*}
\big(\widehat u(t),\widehat b(t),\widehat x(t)\big):=(\bar{u}(0),\bar
b(0),x_0)+\int^t_0\big(\vartheta^u(\tau),\vartheta^b(\tau),\vartheta^x(\tau)%
\big)d\tau\;\mbox{ for all }\;t\in[0,T]
\end{equation*}
and get that $(\dot{\widehat u}(t),\dot{\widehat b}(t),\dot{\widehat x}%
(t))=(\vartheta^u(t),\vartheta^b(t),\vartheta^x(t))$ for a.e.\
$t\in[0,T]$. It follows from the construction of $(\widehat
u(t),\widehat b(t),\widehat x(t))$ and the passage to the limit as
$k\to\infty$ in \eqref{ini}--\eqref{ic2} that $\|\widehat u(t)\|=1$
on $[0,T]$, that $\widehat x(T)\in\Omega$, and that $(\widehat
u(t),\widehat
b(t),\widehat x(t))$ belongs to the $\varepsilon$-neighborhood of $(\bar{u}%
(\cdot),\bar b(\cdot),\bar{x}(\cdot))$ in the norm topology of $%
W^{1,2}([0,T],\mathbb{R}^{mn}\times\mathbb{R}^m\times\mathbb{R}^n)$. Let us
now check that the limiting triple $(\widehat u(\cdot),\widehat
b(\cdot),\widehat x(\cdot))$ satisfies the sweeping inclusion
\begin{equation}  \label{sweep}
-\dot{x}(t)\in N_{C_{(u(t),b(t))}}\big(x(t)\big)\;\mbox{ for a.e.
}\;t\in[0,T]
\end{equation}
over the controlled polyhedron. It follows from 
\eqref{disc-sw} due to \eqref{movpoly} and \eqref{F} that
\begin{equation*}
\left\langle\bar{u}^k_i(t_j),\bar{x}^k(t_j)\right\rangle\le\bar b^k_i(t_j)\;
\mbox{ for all }\;i=1,\ldots,m,\;\mbox{ all }\;j=0,\;\nu(k)-1,\;\mbox{ and }
\;k\in\mathbb{N}.
\end{equation*}
Passing there to the limit as $k\to\infty$ ensures the conditions
\begin{equation}  \label{Ck}
\left\langle\widehat u_i(t),\widehat x(t)\right\rangle\le\widehat b_i(t)\;
\mbox{ for all }\;i=1,\ldots,m\;\mbox{ and }\;t\in[0,T],
\end{equation}
i.e., $\widehat x(t)\in C_{(\widehat u(t),\widehat b(t))}$ on $[0,T]$. To
proceed further, we use the construction of $F$ in \eqref{F} and rewrite 
\eqref{disc-sw} along the optimal triple $(\bar{u}^k,\bar b^k,\bar{x}^k)$
for $(P_k)$ as
\begin{equation}  \label{disc-sw1}
-\frac{\bar{x}^k(t_{j+1})-\bar{x}^k(t_j)}{h_{k_j}}\in N_{C_{(\bar{u}%
^k(t_j),\bar b^k(t_j))}}\big(\bar{x}^k(t_j)\big)\quad (j=0,\ldots,\nu(k)-1,\,\,k\in\mathbb{N}).
\end{equation}
Recalling the piecewise linear extensions $(\bar{u}^k(t),\bar b^k(t),\bar{x}^k(t))$ of the discrete triples $(%
\bar{u}^k,\bar b^k,\bar{x}^k)$ and their strong
$W^{1,2}$-convergence to the triple $(\widehat u(t),\widehat
b(t),\widehat x(t))$ satisfying \eqref{Ck} tells us by passing to
the limit in \eqref{disc-sw1} as $k\to\infty$ that the sweeping
inclusion \eqref{sweep} holds for $(\widehat u(t),\widehat
b(t),\widehat x(t))$. The verification of the latter involves the
usage of the aforementioned Mazur theorem and the outer
semicontinuity (closed-graph) property of the convex normal cone
\eqref{nc} with respect to pointwise perturbations of the moving
polyhedron $C_{(u,b)}$ in \eqref{sweep}.

All the above shows that the limiting triple $(\widehat u,\widehat
b,\widehat x)$ is a feasible solution to problem $(P)$ while
satisfying the $\varepsilon$-localization condition \eqref{loc}.
Passing finally to the limit in $(P_k)$ with taking into account the assumed continuity of $
\varphi$ and remembering the value $\sigma>0$ of the chosen limiting
point of the sequence in \eqref{lim-con}, we get that
$\varphi(\widehat x(T))<\varphi(\bar{x}(T))$. This contradicts the
imposed local optimality of $(\bar{u},\bar b)$ in $(P)$ and hence
completes the proof of theorem. \qed
\end{proof}\vspace*{-0.27in}

\section{Optimality conditions via discrete
approximations}\label{sec:optim-disc}\setcounter{equation}{0}\vspace*{-0.05in}

The results of the previous section show that optimal solutions to
the finite-dimensional discrete-time problem $(P_k)$ are
approximating {\em suboptimal} solutions to the original sweeping
control problem $(P)$ of infinite-dimensional dynamic optimization.
Therefore, necessary optimality conditions for solutions to problems
$(P_k)$, when $k\in\N$ is sufficiently large, can be viewed as
(necessary) {\em suboptimality conditions} for the prescribed
intermediate local minimizers of $(P)$. This observation allows us
to justify solving the original sweeping control problem by applying
appropriate numerical techniques based on necessary optimality
conditions for the discrete approximations.

Each discrete-time problem $(P_k)$ can be reduced to a nondynamic
problem of \emph{mathematical programming} in finite-dimensional spaces. As
we see, problems $(P_k)$ contain constraints of special types, the
most challenging of which are given by \emph{increasingly many inclusions} in
\eqref{disc-sw} that come from the sweeping dynamics. The latter
constraints of the \emph{graphical type} require appropriate tools of generalized differentiation to deal with. In particular, Clarke's nonsmooth analysis cannot be apply here, since his normal cone is usually too large for
graphical sets associated with velocity mappings in
\eqref{sweeping} and \eqref{disc-sw}. In fact, the only (known to
us) machinery of generalized differentiation suitable for these
purposes is the one introduced by the third author and then
developed by many researchers; see, e.g., the books
\cite{m06,m18,rw} and the references therein. We now briefly review
what is needed in this paper.

Given a set $\Th\subset\R^n$ locally closed around $\oz\in\Th$, the
(Mordukhovich basic/limiting) {\em normal cone} to $\Th$ at $\oz$ is
defined by
\begin{align}\label{nor_con}
&N(\oz;\Th)=N_\Th(\oz):=\\
&\big\{v\in\R^n\;\big|\;\exists\,z_k\to\oz,\;w_k\in\Pi(z_k;\O),\;\al_k\ge
0\;\mbox{ with }\;\al_k(z_k-w_k)\to v\big\},\notag
\end{align}
where $\Pi(z;\Th):=\{w\in\Th\;|\;\|z-w\|=d(z,\Th)\}$ is the
Euclidean projector of $z\in\R^n$ onto $\Th$. While for convex sets
$\Th$ the normal cone \eqref{nor_con} agrees with the classical one
\eqref{nc}, in general the set of normals \eqref{nor_con} may be
nonconvex even for simple sets as, e.g., the graph of the absolute
value function $|\cdot|$ at $\oz=(0,0)\in\R^2$. Nevertheless, the
normal cone \eqref{nor_con} for sets, as well as the coderivatives
of set-valued mappings and (first-order and second-order)
subdifferentials of extended-real-valued functions generated by
\eqref{nor_con}, enjoy \emph{comprehensive calculus rules} that are based
on \emph{variational and extremal principles} of variational analysis.

Given further a set-valued mapping ${\cal F}\colon\R^n\tto\R^m$ with the
graph $\gph{\cal F}:=\{(x,y)\in\R^n\times\R^m\;|\;y\in{\cal F}(x)\}$
locally closed around $(\ox,\oy)\in\gph{\cal F}$, the {\em
coderivative} of ${\cal F}$ at $(\ox,\oy)$ is defined by
\begin{equation}\label{coderivative}
D^*{\cal F}(\ox,\oy)(u):=\big\{v\in\R^n\;\big|\;(v,-u)\in
N\big((\ox,\oy);\gph{\cal F}\big)\big\},\quad u\in\R^m.
\end{equation}
Given finally an extended-real-valued function
$f\colon\R^n\to\oR:=(-\infty,\infty]$ lower semicontinuous around
$\ox$ with $f(\ox)<\infty$ and the epigraph $\epi
f:=\{(x,\al)\in\R^{n+1}\;|\;\al\ge f(x)\}$, the (first-order) {\em
subdifferential} of $f$ at $\ox$ can be defined geometrically via
the normal cone \eqref{nor_con} as
\begin{equation}\label{sub}
\partial f(\ox):=\big\{v\in\R^n\big|\;(v,-1)\in
N\big((\ox,f(\ox));\epi f\big)\big\},
\end{equation}
while it admits various analytic descriptions that can be found in
the aforementioned books. Observe that the normal cone
\eqref{nor_con} is the subdifferential \eqref{sub} of the indicator
function $\dd_\Th(x)$ of $\Th$, which equals $0$ for $x\in\Th$ and
$\infty$ otherwise. The {\em second-order subdifferential} of $f$ at
$\ox$ relative to $\ox\in\partial f(\ox)$ is defined as the
coderivative of the first-order subdifferential mapping by
\begin{equation}\label{2nd}
\partial^2f(\ox,\ov)(d):=\big(D^*\partial f\big)(\ox,\ov)(d),\quad
d\in\R^n.
\end{equation}
This construction naturally arises in optimal control of sweeping
processes of type \eqref{sweeping}, where the right-hand side is
described by the normal cone mapping. We look for
second-order evaluations of the coderivative in \eqref{2nd} applied
to the normal cone mapping $F$ in \eqref{F} generated by the
control-dependent convex polyhedron $C(u,b)$ in the sweeping process
\eqref{sweeping}. The result needed in this paper follows from
\cite[Theorem~4.3]{mor2}, where it was derived by using calculations
in \cite{mo} and Robinson's theorem of the calmness property
of polyhedral multifunctions \cite{rob}.
To formulate the required result, consider the matrix 
$$
A:=[u_{ij}]\,\,(i=1,\ldots,m;j=1,\ldots,n)
$$
with the vector columns $u_i$ as well
as the transpose matrix $A^T$. As usual, the symbol $^\perp$
indicates the orthogonal complement of a vector in the corresponding
space. Having the controlled polyhedron $C(u,b)$ in \eqref{F}, take
its {\em active indices} at $(u,b,x)$ with $x\in C(u,b)$ denoted by
\begin{equation*}
I(u,b,x):=\big\{i\in\{1,\ldots,m\}\;\big|\;\la u_i,x\ra=b_i\big\}.
\end{equation*}
The {\em positive linear independence constraint
qualification} (PLICQ) at $(u,b,x)$ is
\begin{equation}\label{PLICQ}
\bigg[\sum_{i\in I(x,u,b)}\al_iu_i=0,\,\al_i\ge
0\bigg]\Longrightarrow\big[\al_i=0\;\;\mbox{for all}\;\;i\in
I(x,u,b)\big].
\end{equation}
This condition is significantly weaker than the classical {\em linear independence constraint qualification} (LICQ), which corresponds to \eqref{PLICQ} with
$\al_i\in\R$ while not being used in this paper. Considering the
moving polyhedron as in \eqref{movpoly}, it is not hard to
check that our basic uniform Slater condition from \eqref{unifslater}
is equivalent to the fulfillment of PLICQ along the feasible triple
$(x(t),u(t),b(t))$ for all $t\in[0,T]$; see \cite{mor1} for more
discussions on this topic.

Given $x\in C(u,b)$ and $v\in N(x;C(u,b))$, define the sets
\begin{equation*}
Q(p):=\left\{\begin{array}{ll}
q_i=0\;\mbox{ for all }\;i\;\mbox{ with either }\;\la u_i,x\ra<b_i\;\mbox{ or }\;p_i=0,\;\mbox{ or}\;\la u_i,y\ra<0,\\
q_i\ge 0\;\mbox{ for all }\;i\;\mbox{ such that }\;\la
u_i,x\ra=b_i,\;p_i=0,\;\mbox{ and }\;\la u_i,y\ra>0,
\end{array}
\right.
\end{equation*}
\begin{equation*}
P(y):=\big\{p\in N_{\R^m_-}(Ax-b)\;\big|\;A^Tp=v\big\}\;\mbox{ for
}\;y\in\bigcap_{\{i\;|\;p_i>0\}}u_i^\perp,
\end{equation*}
where the normal cone to the nonpositive orthant $\R^m_-$ is easy to
compute.

Now we are ready to present the required evaluation of the
coderivative of the normal cone mapping $F(x,u,b)$ generated by the
controlled polyhedron in \eqref{F}. The following lemma is a slight
modification of \cite[Theorem~4.3]{mor2}.\vspace*{-0.1in}

\begin{lemma}\label{cod-eval} Taking $F$ and $C(u,b)$ from
\eqref{F}, suppose that the active vector columns $\{u_i\;|\;i\in
I(u,b,x)\}$ are positively linearly independent for any $(u,b,x)$
with $x\in C(u,b)$. Then for all such $(u,b,x)$, all $v\in
N(x;C(u,b))$, and all $y\in\cap_{\{i\;|\;p_i>0\}}u_i^\perp$ we have
the coderivative upper estimate
\begin{equation}\label{cod_inclusion}
D^*F(u,b,x,v)(y)\subset\bigcup\limits_{\begin{subarray}{l}p\in
P(y)\\q\in Q(p)
\end{subarray}}\left\{\left(\begin{array}{c}
A^Tq\\\hline
p_{1}y+q_1x\\
\vdots\\
p_my+q_mx\\\hline-q
\end{array}
\right)\right\}.
\end{equation}
\end{lemma}

Note that imposing the LICQ condition instead of PLICQ ensures that
the set $P(y)$ is a singleton and that the inclusion in
\eqref{cod_inclusion} holds as equality; see \cite[Theorem~4.3]{mor2}. However, for the purpose of this paper it is sufficient to have the inclusion in \eqref{cod_inclusion} under the less restrictive PLICQ. 

To proceed further, we need one more auxiliary result giving us
necessary optimality conditions for a finite-dimensional nondynamic
problem of {\em mathematical programming} with finitely many
equality, inequality and inclusion (geometric) constraints. The next
lemma is obtained by combining the necessary optimality conditions
from \cite[Theorem~6.5]{m18} for mathematical programs containing
one geometric constraint and the intersection rule for limiting
normals taken from \cite[Corollary~2.17]{m18}. Arguing in this way,
we can derive necessary optimality conditions for mathematical
programs described by lower semicontinuous cost and inequality
constraint functions as well as continuous functions describing
equality constraints. However, we confine ourselves to considering
problems with just locally Lipschitzian functions for cost and
inequality constraints and smooth functions for equality
constraints, since only such functional constraints appear in
mathematical programs to which we reduce the discrete-time sweeping
control problems $(P_k)$.\vspace*{-0.1in}

\begin{lemma}\label{math-prog} Consider the following problem of
mathematical programming:
\begin{equation} \tag{MP}
\left\{\begin{array}{ll}\mbox{minimize }\;f_0(z)\;\mbox{ as
}\;z\in\R^d\;\mbox{ subject
to}\\
f_i(z)\le 0\;\mbox{ for }\;i=1,\ldots,s,\\
g_j(z)=0\;\mbox{ for }\;j=0,\ldots,r,\\
z\in\Th_j\;\mbox{ for }\;j=0,\ldots,l,
\end{array}\right.
\end{equation}
where all the functions $f_i$ and $g_j$ are real-valued. Given a
local minimizer $\oz$ to $(MP)$, assume that the functions $f_i$ are
locally Lipschitzian around $\oz$ for $i=0,\ldots,s$, the functions
$g_j$ are continuously differentiable around this point for
$j=0,\ldots,r$, and the sets $\Th_j$ are locally closed around $\oz$
for all $j=0,\ldots,l$. Then there exist nonnegative numbers
$\lm_0,\ldots,\lm_{s}$, real numbers $\mu_0,\ldots,\mu_r$, and
vectors $z^*_j\in\R^d$ for $j=0,\ldots,l$, not equal to zero
simultaneously, such that
\begin{equation*}
\lm_i f_i(\oz)=0\;\mbox{ for }\;i=1,\ldots,s,
\end{equation*}
\begin{equation*}
z^*_j\in N(\oz;\Th_j)\;\mbox{ for }\;j=0,\ldots,l,
\end{equation*}
\begin{equation*}
-\sum_{j=0}^l z^*_j\in\sum_{i=0}^s\lm_i\partial
f_i(\oz)+\sum_{j=0}^{r}\mu_i\nabla g_j(\oz),
\end{equation*}
\end{lemma}

Having Lemma~\ref{cod-eval} and Lemma~\ref{math-prog} in hand, we
are now in a position to establish necessary conditions for optimal
solutions to problems $(P_k)$ from \eqref{disc-cost}--\eqref{ic2}
whenever the approximation number $k\in\N$ is sufficiently large.
The obtained relationships involve the given intermediate local
minimizer for the sweeping optimal control problem $(P)$ and thus
present necessary suboptimality conditions for the original
continuous-time problem due to Theorem~\ref{ilm-conver}. For any
$x\in\R^n$, $y=(y_1,\ldots,y_m)\in\R^{nm}$ with $y_i\in\R^n$
$(i=1,\ldots,m)$, and $\al=(\al_1,\ldots,\al_m)\in\R^m$ we use the
symbols
\begin{equation*}
{\rm rep}_m(x):=(x,\ldots,x)\in\R^{nm}\;\mbox{ and
}\;[\alpha,y]:=(\alpha_1y_1,\ldots,\alpha_my_m)\in\R^{nm}.
\end{equation*}

\begin{theorem}\label{nc-disc} Let $(\ou,\ob)$ be a canonical
intermediate local minimizer of $(P)$ generated the trajectory
$\ox=\ox(\cdot)$ of the controlled polyhedral sweeping process
\eqref{sweeping} such that the cost function $\ph$ is locally
Lipschitzian around $\ox(T)$. Fix an optimal triple
$(\ou^k,\ob^k,\ox^k)$  in problem $(P_k)$ with the components
\begin{equation*}
(\ou^k,\ob^k,\ox^k):=(\ou^k_0,\ou^k_1,\ldots,\ou^k_{\nu(k)},\ob^k_0,\ob^k_1,\ldots,\ob^k_{\nu(k)},\ox^k_0,\ox^k_1,\ldots,\ox^k_{\nu(k)})
\end{equation*}
and choose $k\in\N$ to be sufficiently large. Denote the quantities
\begin{align*}
\th_{j}^{uk}&:=\int_{t_j^k}^{t_{j+1}^k}\(\frac{\ou_{j+1}^k-\ou_j^k}{h^k_j}-\dot{\ou}(t)\)dt,\quad
\th_{j}^{bk}:=\int_{t_j^k}^{t_{j+1}^k}\(\frac{\ob_{j+1}^k-\ob_j^k}{h^k_j}-\dot{\ob}(t)\)dt,\\
\th_{j}^{xk}&:=\int_{t_j^k}^{t_{j+1}^k}\(\frac{\ox_{j+1}^k-\ox_j^k}{h^k_j}-\dot{\ox}(t)\)dt
\end{align*}
and define the set $\O_k:=\O+\xi_kI\!\!B$, where $\xi_k$ is taken
from the construction of problem $(P_k)$. Then there exist a
multiplier $\lm^k\ge 0$, an adjoint triple
$p_j^k=(p_{j}^{xk},p_{j}^{ak},p_{j}^{bk})\in\R^{n+mn+m}$
$(j=0,\ldots,\nu(k))$, as well as vectors
$\eta^k=(\eta_0^k,\ldots,\eta^k_{\nu(k})$ $\in\R^{m(\nu(k)+1)}_+$,
$\al^{1k}=\(\al_{0}^{1k},\ldots,\al_{\nu(k)}^{1k}\)\in\R^{m(\nu(k)+1)}_+$,
$\al^{2k}=(\al_{0}^{2k},\ldots,\al_{\nu(k)}^{2k})\in\R^{m(\nu(k)+1)}_{+}$,
and $\gg^k=(\gg_0^k,\ldots,\gg^k_{\nu(k)-1})\in\R^{m\nu(k)}$ such
that
\begin{equation}\label{ntc0}
\lm^k+\n\al^{1k}-\al^{2k}\en+\disp\n\eta_{\nu(k)}^k\en+\sum_{j=0}^{\nu(k)-1}\n
p_{j}^{xk}\en+\n p_{0}^{ak}\en+\n p_{0}^{bk}\en\ne 0,
\end{equation}
\begin{equation}\label{ntc1}
\lm^k+\n\al^{1k}-\al^{2k}\en+\n\gg^k\en+\n p^{ak}_{\nu(k)}\en+\n
p^{bk}_{\nu(k)}\en\ne 0,
\end{equation}
and we have the following conditions:\\[1ex]
{\sc $\bullet$ dynamic relationships}, which are satisfied for all indices
$j=0,\ldots,\nu(k)-1$ and $i=1,\ldots,m:$
\begin{equation}\label{87}
-\frac{\ox_{j+1}^k-\ox_j^k}{h^k_j}
=\sum_{i=1}^{m}\eta_{ij}^k\ou_{ij}^k,
\end{equation}
\begin{equation}\label{cona}
\begin{array}{ll}
&\dfrac{p_{j+1}^{uk}-p_{j}^{uk}}{h^k_j}-\dfrac{2}{h^k_j}\[\al_{j}^{1k}-\al_{j}^{2k},\ou_{j}^k\]\\
&=\[\gg_{j}^k,\rep_m(\ox_j^k)\]+\[\eta_j^k,\rep_m\(-\dfrac{1}{h^k_j}\lm^k\th_{j}^{xk}-\lm^k+p_{j+1}^{xk}\)\],
\end{array}
\end{equation}
\begin{equation}\label{conb}
\frac{p_{j+1}^{bk}-p_{j}^{bk}}{h^k_j}=-\gg_{j}^k,
\end{equation}
\begin{equation}\label{conx}
\frac{p_{j+1}^{xk}-p_{j}^{xk}}{h^k_j}=\sum_{i=1}^{m}\gg_{ij}^k\ou_{ij}^k,
\end{equation}
where the components of the vectors $\gg^k_j$ are such that
\eq\label{congg1}
\begin{cases}
\gg^k_{ij}=0\;\;\mbox{if}\;\;\la\ou^k_{ij},\ox^k_j\ra<\ob^k_{ij}\;\;\mbox{or}\;\;\eta^k_{ij}=0,\;\Big\la\ou_{ij}^k,-\dfrac{1}{h^k_j}\lm^k\th_{j}^{xk}+p_{j+1}^{xk}\Big\ra<0,\\
\gg^k_{ij}\ge 0\;\;\mbox{if}\;\;\la\ou^k_{ij},\ox^k_j\ra=\ob^k_{ij},\;\eta^k_{ij}=0,\;\Big\la\ou_{ij}^k,-\dfrac{1}{h^k_j}\lm^k\th_{j}^{xk}+p_{j+1}^{xk}\Big\ra>0,\\
\gg^k_{ij}\in\R
\;\;\mbox{if}\;\;\eta^k_{ij}>0,\;\Big\la\ou_{ij}^k,-\dfrac{1}{h^k_j}\lm^k\th_{j}^{xk}+p_{j+1}^{xk}\Big\ra=0.
\end{cases}
\eeq {\sc $\bullet$ complementary slackness conditions}:
\begin{equation}\label{71l1}
\al_{ij}^{1k}\(\n u_{ij}^k\en-(1+\dd_k)\)=0\quad (i=1,\ldots,m,\;\;j=0,\ldots,\nu(k)),
\end{equation}
\begin{equation}\label{71l2}
\al_{ij}^{2k}\(\n u_{ij}^k\en-(1-\dd_k)\)=0\quad
(i=1,\ldots,m,\;\;j=0,\ldots,\nu(k)),
\end{equation}
\begin{equation}\label{eta}
\[\la u_{ij}^k,\ox_j^k\ra<\ob_{ij}^k\]\sr\eta_{ij}^k=0\quad (i=1,\ldots,m,\;\;j=0,\ldots,\nu(k)-1),
\end{equation}
\begin{equation}\label{eta1}
\[\la\ou_{i\nu(k)}^k,\ox_{\nu(k)}^k\ra<\ob_{i\nu(k)}^k\]\sr\eta_{i\nu(k)}^k=0\quad (i=1,\ldots,m,\;\;j=0,\ldots,\nu(k)-1),
\end{equation}
\begin{equation}\label{96}
\eta_{ij}^k>0\sr\[\Big\la\ou_{ij}^k,-\frac{1}{h^k_j}\lm^k\th_{j}^{xk}+p_{j+1}^{xk}\Big\ra=0\]\,\,(i=1,\ldots,m,\;\;j=0,\ldots,\nu(k)-1).
\end{equation}
{\sc $\bullet$ transversality relationships} at the right end of the
trajectory:
\begin{equation}\label{nmutx}
-p_{\nu(k)}^{xk}\in\lm^k\partial\ph(\ox_{\nu(k)}^k)+N\big(\ox^k_{\nu(k)};\O_k)+\sum_{i=1}^{m}\eta_{i\nu(k)}^k\ou_{i\nu(k)}^k,
\end{equation}
\begin{equation}\label{nmuta}
p_{\nu(k)}^{uk}=-2\[\al_{\nu(k)}^{1k}-\al_{\nu(k)}^{2k},\ou_{\nu(k)}^k\]-\[\eta_{\nu(k)}^k,\rep_m(\ox_{\nu(k)}^k)\],
\end{equation}
\begin{equation}\label{nmutb}
p_{i\nu(k)}^{bk}=\eta^k_{i\nu(k)}\ge 0,\;\la
\ou^{k}_{i\nu(k)},\ox^k_{\nu(k)}\ra<\ob^k_{i\nu(k)}\sr
p_{i\nu(k)}^{bk}=0\quad (i=1,\ldots,m).
\end{equation}
\end{theorem}\vspace*{-0.15in}
\begin{proof} To reduce problem $(P_k)$ from \eqref{disc-cost}--\eqref{ic2} for each fixed $k\in\N$ to a
mathematical program of type $(MP)$ formulated in
Lemma~\ref{math-prog}, we form the multidimensional vector
\begin{align*}
z:=\left(u_0^k,\ldots,u_{\nu(k)}^k,b_0^k,\ldots,b_{\nu(k)}^k,x_0^k,\ldots,x_{\nu(k)}^k,v_0^k,\ldots,v_{\nu(k)-1}^k,\right.&\\ \left. w_0^k,\ldots,w_{\nu(k)-1}^k,y_0^k,\ldots,y_{\nu(k)-1}^k\right)&
\end{align*}
and consider the problem of minimizing the cost function
\begin{equation*}
f_0(z):=\vph(x^k_{\nu(k)})+\frac{1}{2}\sum_{j=0}^{\nu(k)-1}\int_{t_j^k}^{t_{j+1}^k}\big\|\big(v_j^k-\dot{\ou}(t),w_j^k-\dot{\ob}(t),y_j^k-\dot{\ox}(t)
\big)\big\|^2dt
\end{equation*}
subject to the five groups of inequality constraints
\begin{equation*}
f_1(z):=\sum_{j=0}^{\nu(k)-1}\int_{t^k_j}^{t^k_{j+1}}\n\(u^k_j,b^k_j,x^k_j\)-\(\ou(t),\ob(t),\ox(t)\)\en^2dt-\dfrac{\ve}{2}\le
0,
\end{equation*}
\begin{equation*}
f_2(z):=\sum_{j=0}^{\nu(k)-1}\int_{t_j^k}^{t_{j+1}^k}\n\big(v_j^k,w_j^k,y_j^k\big)-\big(\dot{\ou}(t),\dot{\ob}(t),\dot{\ox}(t)\big)\en^2dt-\frac{\ve}{2}\le
0,
\end{equation*}
\begin{equation*}
f_{ij}(z):=\|u_{ij}^k\|^2-(1+\dd_k)^2\le 0\;\mbox{ for
}\;i=1,\ldots,m,\;\;j=0,\ldots,\nu(k),
\end{equation*}
\begin{equation*}
\tilde f_{ij}(z):=(1-\dd_k)^2-\|u_{ij}^k\|^2\le 0,\;\mbox{ for
}\;i=1,\ldots,m,\;\;j=0,\ldots,\nu(k),
\end{equation*}
\begin{equation*}
\hat f_i(z):=\big\la
u_{i\nu(k)}^k,x_{\nu(k)}^k\big\ra-b_{i\nu(k)}^k\le 0\;\mbox{ for
}\;i=1,\ldots,m,
\end{equation*}
the three groups of equality constraints
\begin{equation*}
g^u_j(z):=u_{j+1}^k-u_j^k-h^k_jv_j^k=0\;\;\mbox{ for
}\;j=0,\ldots,\nu(k)-1,
\end{equation*}
\begin{equation*}
g^b_j(z):=b_{j+1}^k-b_j^k-h^k_jw_j^k=0\;\mbox{ for
}\;j=0,\ldots,\nu(k)-1,
\end{equation*}
\begin{equation*}
g^x_j(z):=x_{j+1}^k-x_j^k-h^k_jy_j^k=0,\;\mbox{ for
}\;j=0,\ldots,\nu(k)-1,
\end{equation*}
and the two groups of inclusion constraints
\begin{equation*}
z\in\Th_j:=\big\{z\;\big|-y_j^k\in
F(u_j^k,b_j^k,x_j^k)\big\}\;\mbox{ for }\;j=0,\ldots,\nu(k)-1,
\end{equation*}
\begin{equation*}
z\in\Th_{\nu(k)}:=\big\{z\;\big|\;(u^k_0,b^k_0,x^k_0)\;\mbox{ are
fixed, }\;x^k_{\nu(k)}\in\O_k\big\},
\end{equation*}
where those for $j=0,\ldots,\nu(k)-1$ incorporate the constraints
$x^k_j\in C(u^k_j,b^k_j)$ for such $j$ due to the construction of
$F$ in \eqref{F}.

As we see, the formulated nondynamic equivalent of problem $(P_k)$
is written in the mathematical programming form $(MP)$ as in
Lemma~\ref{math-prog} with the fulfillment all the assumptions
imposed in the lemma. Thus we can readily apply the conclusions of
the lemma by taking into account the particular structure of the
functions and sets in the formulated equivalent of $(P_k)$.
Employing now the necessary optimality conditions of
Lemma~\ref{math-prog} to the optimal solution
\begin{align*}
\oz:=\oz^k=\(\ou_0^k,\ldots,\ou_{\nu(k)}^k,\ob_0^k,\ldots,\ob_{\nu(k)}^k,\ox_0^k,\ldots,\ox_{\nu(k)}^k,\ov_0^k,\ldots,\ov_{\nu(k)-1}^k,\right.&\\ \left.\ow_0^k,\ldots,\ow_{\nu(k)-1}^k,
\oy_0^k,\ldots,\oy_{\nu(k)-1}^k\)&
\end{align*}
of problem $(MP)\equiv(P_k)$, observe by Theorem~\ref{ilm-conver}
that the inequality constraints defined by the functions $f_1$ and
$f_2$ above are {\em inactive} at $\oz$ for sufficiently large $k$,
and thus the corresponding multipliers will not appear in optimality
conditions. Taking this into account, we find by
Lemma~\ref{math-prog} multipliers $\lm^k\ge 0$,
$(\bb^k_1,\ldots,\bb^k_m)\in\R^m_+$,
$p^k_j=(p^{uk}_j,p^{bk}_j,p^{xk}_j)\in\R^{mn+n+m}$ for
$j=1,\ldots,\nu(k)$, as well as vectors
\begin{align*}
z^*_j:=\(u^*_{0j},\ldots,u^*_{\nu(k)j},b^*_{0j},\ldots,b^*_{\nu(k)j},x^*_{0j},\ldots,x^*_{\nu(k)j},v^*_{0j},\ldots,v^*_{(\nu(k)-1)j},\right.&\\ \left. w^*_{0j},\ldots,w^*_{(\nu(k)-1)j},y^*_{0j},\ldots,y^*_{(\nu(k)-1)j}\)&
\end{align*} 
for $j=0,\ldots,\nu(k)$,
$\al^{1k}=(\al_{0}^{1k},\ldots,\al_{\nu(k)}^{1k})\in\R^{\nu(k)+1}_+$,
$\al^{2k}=(\al_{0}^{2k},\ldots,\al_{\nu(k)}^{2k})\in\R^{\nu(k)+1}_{-}$
such that the complementary slackness conditions \eqref{71l1},
\eqref{71l2}, and
\begin{equation}\label{71+}
\bb^k_i\big(\big\la
\ou_{i\nu(k)}^k,\ox_{\nu(k)}^k\big\ra-\ob_{i\nu(k)}^k\big)=0\;\mbox{
for }\;i=1,\ldots,m
\end{equation}
hold together with the normal cone inclusions
\begin{equation}\label{nor-inc}
z^*_j\in N(\oz;\Th_j)\;\mbox{ for }\;j=0,\ldots,\nu(k)
\end{equation}
and the generalized Lagrangian condition
\begin{equation}\label{70}
\begin{array}{ll}
-\disp\sum_{j=0}^{\nu(k)}z^*_j&\in\lm^k\partial
f_0(\oz)+\disp\sum_{i=1}^{m}\bb^k_j\nabla\hat
f_i(\oz)+\disp\sum_{j=0}^{\nu(k)-1}\nabla g_j(\oz)^T p_{j+1}^{k}\\
&+\disp\sum_{j=0}^{\nu(k)}\sum^m_{i=1}\bigg[\al_{ij}^{1k}\nabla
f_{ij}(\oz)+\al_{ij}^{2k}\nabla\tilde f_{ij}(\oz)\bigg],
\end{array}
\end{equation}
where $g_j=(g^u_j,g^b_j,g^x_j)$, and where the dual elements $\lm^k$, $\bb^k_i$,
$p^k_j$, $z^*_j$, $\al^{1k}$, and $\al^{2k}$ are not all zero
simultaneously.

Looking at the graphical structure of the geometric constraints
$z\in\Th_j$ for $j=0,\ldots,\nu(k)-1$, we readily deduce from
\eqref{nor-inc} that
\begin{equation*}
\begin{array}{ll}
(u^*_{jj},b^*_{jj},x^*_{jj},-y^*_{jj})\in
N\Big(\Big(\ou_j^k,\ob_j^k,\ox_j^k,-\dfrac{\ox_{j+1}^k-\ox_j^k}{h^k_j}
\Big);\gph F\Big)\,\,(j=0,\ldots,\nu(k)-1)
\end{array}
\end{equation*}
with all the other components of $z^*_j$ equal to zero for these
indices $j$. It follows from the coderivative definition
\eqref{coderivative} that the obtained normal cone inclusion can be
equivalently written as
\begin{equation}\label{cod-inc}
(u^*_{jj},b^*_{jj},x^*_{jj})\in
D^*F\Big(\ou_j^k,\ob_j^k,\ox_j^k,-\dfrac{\ox_{j+1}^k-\ox_j^k}{h^k_j}
\Big)(y^*_{jj})\;\mbox{ for }\;j=0,\ldots,\nu(k)-1.
\end{equation}
Since the mapping $F$ is given in the particular form \eqref{F}, we
are able to use the coderivative evaluation in \eqref{cod-inc}
provided the fulfillment of PLICQ \eqref{PLICQ} along the discrete
optimal solutions for all $k$ sufficiently large. As discussed
above, the assumed uniform Slater condition \eqref{unifslater} for
the given canonical intermediate local minimizer $(\ou,\ob)$ of
$(P)$ yields PLICQ at $(\ou,\ob,\ox)$. Since the latter condition is
{\em robust} with respect to perturbations of the initial triple and
since the discrete optimal solutions strongly converge to
$(\ou(\cdot),\ob(\cdot),\ox(\cdot))$ by Theorem~\ref{ilm-conver}, we
are in a position to use Lemma~\ref{cod-eval} in the coderivative
inclusion \eqref{cod-inc}. Prior to this, let us calculate the other
terms in the generalized Lagrangian condition \eqref{70}.

First observe that the summation term in the cost function is
smooth. Therefore, the usage of the subdifferential sum rule from
\cite[Proposition~1.30(ii)]{m18} gives the precise calculation
\begin{equation*}
\partial
f_0(\oz)=\partial\ph(\ox^k_{\nu(k)})+\sum_{j=0}^{\nu(k)-1}\big(0,\ldots,0,\th^{uk}_j,\th^{bk}_j,\th^{xk}_j\big)
\end{equation*}
where zeros stands for the all components of $\oz$ till $\ov^k_j$,
and where $\th^{uk}_j,\th^{bk}_j,\th^{xk}_j$ are defined in the
formulation of the theorem. Further, with the usage of our notation
presented before the formulation of this theorem, we easily get
\begin{equation*}
\sum_{i=1}^{m}\bb_i^k\nabla \Hat f_i(\oz)=\(\sum_{i=1}^{m}
\bb_i^k\ou_{ik}^k,\[\bb^k,\rep_m(\ox_{\nu(k)}^k)\],-\bb^k\),
\end{equation*}
\begin{equation*}
\(\sum_{j=0}^{\nu(k)-1}\nabla
g_j(\oz)^Tp_{j+1}^k\)_{(u_j,b_j,x_j)}=\left\{\begin{array}{llll}
-p_{1}^{k}&\textrm{ if }& j=0\\[1ex]
p_{j}^{k}-p_{j+1}^{k}&\textrm{ if }& j=1,\ldots,\nu(k)-1\\[1ex]
p_{\nu(k)}^{k}&\textrm{ if }& j=\nu(k)
\end{array}
\right.,
\end{equation*}
\begin{align*}
\(\sum_{j=0}^{\nu(k)-1}\nabla
g_j(\oz)^Tp_{j+1}^k\)_{(v_j,w_j,y_j)}=\(
-h^k_0p_{1}^{uk},-h^k_1p_{2}^{uk},\ldots,-h^k_{\nu(k)-1}p_{\nu(k)}^{uk},\right.&\\
\left. -h^k_0p_{1}^{bk},-h^k_1p_{2}^{bk},\ldots,
-h^k_{\nu(k)-1}p_{\nu(k)}^{bk},-h^k_0p_{1}^{xk},-h^k_1p_{2}^{xk},\ldots,
-h^k_{\nu(k)-1}p_{\nu(k)}^{xk}\),&
\end{align*}
\begin{align*}
\sum_{j=0}^{\nu(k)}\sum^m_{i=1}\al_{ij}^{1k}\nabla
f_{ij}(\oz)=2\[\al_{j}^{1k},\ou_{j}^k\],\;
\;\sum_{j=0}^{\nu(k)}\sum^m_{i=1}\al_{ij}^{2k}\nabla\tilde
f_{ij}(\oz)=-2\[\al_{j}^{2k},\ou_{j}^k\]&\\
(j=0,\ldots,\nu(k)).&
\end{align*}
To proceed with \eqref{70}, it remains to express the dual element
$z^*_{\nu(k)}\in N(\oz;\Th_{\nu(k)})$ in \eqref{nor-inc}
corresponding the last geometric constraint
$\oz_{\nu(k)}\in\Th_{\nu(k)}$ in terms of the data of $(P_k)$. We
directly conclude from the structure of $\Th_{\nu(k)}$ that the
components of $z^*_{\nu(k)}$ corresponding to $(u_0^k,b_0^k,x_0^k)$
are free (i.e., just belong to $\R^{mn}\times\R^m\times\R^n$), that
$x^*_{\nu(k)\nu(k)}\in N(\ox^k_{\nu(k)};\O_k)$, and that all the other
components are equal to zero. The fulfillment of PLICQ along
$(\ou^k,\ob^k,\ox^k)\}$ for all $k$ sufficiently large allows us to
find unique vectors $\eta_j^k\in\R^m_+$ such that
\begin{equation*}\label{h:5.39}
\sum_{i=1}^m\eta_{ij}^k\ou_{ij}^k=-\frac{\ox_{j+1}^k-\ox_j^k}{h^k_j}\;\mbox{
for }\;j=0,\ldots,\nu(k)-1.
\end{equation*}
For the last index $j=\nu(k)$, we  put
$\eta^k_{\nu(k)}:=\bb^k\in\R^m_+$. Substituting all the above into
the Lagrangian inclusion \eqref{70} with taking into account the
coderivative upper estimate from Lemma~\ref{cod-eval} gives us the
claimed necessary optimality conditions \eqref{87}--\eqref{nmutb}.
Finally, the nontriviality conditions in \eqref{ntc0} and
\eqref{ntc1} follows directly from \eqref{87}--\eqref{nmutb} and the
nontriviality of the dual elements in Lemma~\ref{math-prog} for the
mathematical program $(MP)$ equivalent to $(P_k)$. Therefore, we
complete the proof of the theorem. \qed
\end{proof}\vspace*{-0.3in}

\section*{Acknowledgment}\vspace*{-0.05in}

The first author acknowledges support by the FMJH Program Gaspard
Monge in optimization and operations research including support to
this program by EDF as well as by the Deutsche
Forschungsgemeinschaft for their support of project B04 within the CRC/Transregio 154.
The work of the second author was partly supported by the EIPHI
Graduate School (contract ANR-17-EURE-0002). Research of the third
author was partly supported by the US National Science Foundation
under grants DMS-1512846 and DMS-1808978, by the US Air Force Office
of Scientific Research grant \#15RT04, and by the Australian
Research Council under grant DP-190100555.

\end{document}